%% file: ijcai22_arxiv.tex
\documentclass{article}
\usepackage{ijcai22}

\usepackage{times}
\usepackage{soul}
\usepackage{url}
\usepackage[hidelinks]{hyperref}
\usepackage[utf8]{inputenc}
\usepackage[small]{caption}
\usepackage{graphicx}
\usepackage{amsmath}
\usepackage{amsthm}
\usepackage{booktabs}
\usepackage{algorithm}
\usepackage{algorithmic}
\usepackage{amsfonts}    
\usepackage{multirow}

\usepackage{pifont}
\usepackage{nccmath}

\usepackage[utf8]{inputenc} % allow utf-8 input
\usepackage[T1]{fontenc}    % use 8-bit T1 fonts
\usepackage{hyperref}       % hyperlinks
\usepackage{amsthm}
\usepackage{amsmath}
\usepackage{url}            % simple URL typesetting
\usepackage{booktabs}       % professional-quality tables
\usepackage{amsfonts}       % blackboard math symbols
\usepackage{nicefrac}       % compact symbols for 1/2, etc.
\usepackage{microtype}      % microtypography
\usepackage{xcolor}         % colors
\usepackage{array}
\usepackage{subfigure}
\usepackage{graphicx}
\usepackage{tikz}
\usepackage{pgf, tikz}
\usetikzlibrary{arrows, automata}
\usepackage{comment}
\usepackage[skins]{tcolorbox}
\usepackage{algorithm}
\usepackage{algorithmic}

\newtheorem{example}{Example}
\newtheorem{theorem}{Theorem}

\newtheorem{lemma}[theorem]{Lemma}
\newtheorem{proposition}[theorem]{Proposition}
\newtheorem{corollary}[theorem]{Corollary}
\newtheorem{assumption}{Assumption}

\usepackage{xcolor}
\usepackage{pgfplots}
\usepackage{mathtools}

\usepackage{xr}
\newcommand{\cmark}{\ding{51}}%
\newcommand{\xmark}{\ding{55}}%

\newtheorem{innercustomgeneric}{\customgenericname}
\providecommand{\customgenericname}{}
\newcommand{\newcustomtheorem}[2]{%
  \newenvironment{#1}[1]
  {%
   \renewcommand\customgenericname{#2}%
   \renewcommand\theinnercustomgeneric{##1}%
   \innercustomgeneric
  }
  {\endinnercustomgeneric}
}

\newcustomtheorem{customthm}{Theorem}
\newcustomtheorem{customlemma}{Lemma}

\newtheorem{remark}[theorem]{Remark}

\newcommand{\diag}{\mathop{\mathbf{diag}}}
\newcommand{\abs}[1]{\ensuremath{\left|#1\right|}}
\newcommand{\norm}[2][]{\ensuremath{\left\Vert #2 \right\Vert}}

\renewcommand{\vec}[1]{\mathbf{#1}}
\newcommand{\vect}[1]{\boldsymbol{\mathbf{#1}}}
\newcommand{\argmin}{\mathop{\mathbf{argmin}}}

\newcommand{\grad}{\mathrm{grad}}

\newcommand{\Retr}{\mathrm{Retr}}

\newcommand{\Hess}{\mathrm{Hess}}
\newcommand{\Exp}{\mathrm{Exp}}
\newcommand{\Log}{\mathrm{Log}}

     \newcommand{\BN}{{\mathbb {N}}}
     
     \newcommand{\BR}{{\mathbb {R}}}

\title{Accelerated Multiplicative Weights Update Avoids Saddle Points almost always}

\author{
Yi Feng$^1$
\and
Ioannis Panageas$^2$\and
Xiao Wang$^{1}$
\affiliations
$^1$Shanghai University of Finance and Economics\\
$^2$University of California, Irvine\\
\emails
2021310186@live.sufe.edu.cn,
ipanagea@ics.uci.edu,
wangxiao@sufe.edu.cn
}

\begin{document}
\maketitle

\begin{abstract}
We consider non-convex optimization problems with constraint that is a product of simplices. A commonly used algorithm in solving this type of problem is the Multiplicative Weights Update (MWU), an algorithm that is widely used in game theory, machine learning and multi-agent systems. Despite it has been known that MWU avoids saddle points \cite{ICMLPPW19}, there is a question that remains unaddressed:  ``Is there an accelerated version of MWU that avoids saddle points provably?'' In this paper we provide a positive answer to above question. We provide an accelerated MWU based on Riemannian Accelerated Gradient Descent in \cite{ZS18}, and prove that the Riemannian Accelerated Gradient Descent, thus the accelerated MWU, almost always avoid saddle points.
%We consider non-convex optimization problem over product manifolds, especially with application in Multiplicative Weights Update over product of simplices. We firstly prove that the Riemannian Accelerated Gradient Descent proposed by \cite{ZS18} can be used in product manifold and avoids saddle points with random initialization. As an application, we provide an accelerated version of Multiplicative Weights Update algorithm which is widely used in game theory, constrained optimization, multi-agent systems and machine learning, and our result on the saddle avoidance of Riemannian Accelerated Gradient Descent implies that the Accelerated Multiplicative Weights Update avoid saddle points.

 %Multiplicative Weights Update (MWU) has been studied extensively due to many applications in constrained optimization, game theory, and machine learning. Despite the popularity of the algorithm, there is a question that remain unaddressed:  ``is there an accelerated version of MWU that converges to second-order stationarity provably ?'' In this paper we provide a positive answer to above question. We focuse on a spherical framework for MWU, i.e., MWU is formulated as a manifold gradient descent on sphere. As two important applications of this approach, we provide an accelerated version of MWU, and we prove it converges to a second-order stationary points almost always.  
 \end{abstract}

%Multiplicative Weights Update (MWU) has beed studied extensively due to many applications in constrained optimization, game theory, and machine learning. Despite the popularity of the algorithm, there are two questions remain unaddressed: (1) whether there is an accelerated version of MWU, and (2) whether there is a perturbed version of MWU so that the efficiency in escaping saddle points can be discussed. This paper focuses on a spherical framework for MWU, i.e., MWU is formulated as a manifold gradient descent. We provide two important applications of this approach: acceleration and non-asymptotic convergence. As for the accelerated MWU, we prove that it converges to second-order stationary points almost always, and for the perturbed MWU, we show that the efficiency in escaping saddle points fits into the framework of perturbed Riemannian gradient descent.

\input{intronew}

\input{mx}

\newpage
\input{expnew}

\section{Conclusion}
In this paper we study an Accelerated Multiplicative Weights Update from the Riemannian geometric viewpoint. We prove that with locally geometric convexity in local minima, RAGD avoids saddle points with random initialization, which implies that A-MWU avoids saddle points. This indicates that A-MWU converges to local minima provided RAGD converges. Our experiments have verified the efficiency of A-MWU in escaping saddle points. 
%In this paper we derive an accelerated MWU from the view point of geometry : the multiplicative weights update algorithm is equivalent to the gradient descent on sphere, and we also prove that accelerated Riemannian gradient descent will avoid saddle points almost always, thus converge to a local minimal almost always. So far we only focus on the optimization aspect of the MWU algorithm (saddle evasion, convergence rate, acceleration, etc.). MWU is also a popular algorithm in game theory, especially in zero-sum games and min-max optimization problems, while the non-convex optimization problem discussed in this paper is the potential game counterpart. Understanding the behavior of A-MWU in min-max optimization is an important work for the future. For example, MWU has been proven divergent and chaotic in zero-sum games in \cite{BG18,ChGP19}, then what is provable guarantee of A-MWU in these cases? The limitation of our result is that we only have asymptotic guarantee of avoiding saddle points, to obtain a specific rate with a simplified (if possible) version of the algorithm in \cite{Boumal2020} is of importance.

%\newpage
\twocolumn
\bibliographystyle{named}
\bibliography{Bibli}

\onecolumn
\section{Appendix}
\subsection{Multi-agent A-MWU}

In this section, we provide the full Multi-agent Accelerated Multiplicative Weights Update. To the full generality, it is allowed that different agent has different set of parameters. 
\begin{assumption}
Let $\mathcal{M}=\Delta_1\times...\times\Delta_n$, where 
\[
\Delta_i = \left\{\vec{x}_i\in\mathbb{R}^{d_i}:\sum_{s=1}^{d_i}x_{is}=1,x_{is}\ge 0\right\}.
\]
Let $\{\alpha_t^i\}_{t\in\mathbb{N}}$, $\{\beta_t^i\}_{t\in\mathbb{N}}$ be sequences for each $i\in[n]$, and it holds that
\[
0<c\le\alpha_{t}^i\le\frac{1}{L}, \ \ \ \text{and} \ \ \ \beta_{t}^i>0
\]
for all $i\in[n]$.
\end{assumption}
Recall that we use $\overrightarrow{\vec{x}}=(\vec{x}_1,...,\vec{x}_n)$ to denote the point in the product space $\mathcal{M}=\Delta_1\times...\times\Delta_n$. The algorithm (Multi-agent A-MWU) is give as follows.

 \begin{algorithm}

%\caption{Variable Step RAGD, \cite{ZS18}}
%\label{alg:C}
\begin{algorithmic}
\STATE {input : $\overrightarrow{\vec{x}}_0, \overrightarrow{\vec{v}}_0, 0<c \le\alpha_t^i< \frac{1}{L}, \beta_t^i >0$,$\delta >0$}, for all $i\in[n]$,
\\
\REPEAT
\STATE Compute $s_t^i\in(0,1)$ from the equation $(s_t^i)^2=\alpha_t^i((1-s_t^i)\gamma_t^i+s_t^i\mu)$.
\\
Set $\bar{\gamma}_{t+1}^i=(1-s_t^i)\gamma_t^i+s_t^i\mu$, $\gamma_{t+1}^i=\frac{1}{1+\beta_t^i}\bar{\gamma}_{t+1}^i$
\\
Compute
\begin{fleqn}
\[
S_i=\left(\prod_{k=1}^{d_i}\frac{x_{ik}(t)}{y_{ik}(t)}\right)^{1/d_i}.
\]
\end{fleqn}
Set 
\begin{fleqn}
\[
y_{ik}(t)=\frac{x_{ik}(t)\exp\left(\frac{s_t^i\gamma_t^i}{\gamma_t^i+s_t^i\mu}\ln\left(S_i\frac{v_{ik}(t)}{x_{ik}(t)}\right)\right)}{\sum_{j=1}^{d_i}x_{ij}\exp\left(\frac{s_t^i\gamma_t^i}{\gamma_t^i+s_t^i\mu}\ln\left(S_i\frac{v_{ij}(t)}{x_{ij}(t)}\right)\right)}
\]
\end{fleqn}

Set 
\begin{fleqn}
\[
x_{ik}(t+1)=y_{ik}(t)\frac{1-\alpha_t^i\frac{\partial f}{\partial x_{ik}}(\overrightarrow{\vec{y}}_t)}{1-\alpha_t^i\sum_{j=1}^{d_i}\frac{\partial f}{\partial x_{ij}}(\overrightarrow{\vec{y}}_t)}
\]
\end{fleqn}
Compute
\begin{fleqn}
\[
S_i'=\left(\prod_{k=1}^{d_i}\frac{y_{ik}(t)}{v_{ik}(t)}\right)^{1/d_i} \ \ \text{for all}\ \ i\in[n].
\]
\end{fleqn}
Set 
\begin{fleqn}
\[
u_{ik}=\frac{(1-s_{t^i})\gamma_t^i}{\bar{\gamma}_t^i}\ln\left(S_i'\frac{v_{ik}(t)}{y_{ik}(t)}\right)+y_{ik}(t)\frac{1-\alpha_t^i\frac{\partial f}{\partial x_{ik}}(\overrightarrow{\vec{y}}_t)}{1-\alpha_t^i\sum_{j=1}^{d_i}y_{ij}(t)\frac{\partial f}{\partial x_{ik}}(\overrightarrow{\vec{y}}_t)}-y_{ik}(t)
\]
\end{fleqn}
Set
\begin{fleqn}
\[
v_{ik}(t+1)=\frac{y_{ik}(t)\exp(u_{ik})}{\sum_{j=1}^{d_i}y_{ik}(t)\exp(u_{ik})}
\]
\end{fleqn}

\UNTIL{$\norm{\grad f(\overrightarrow{\vec{y}}_t)} \le \delta$}. %for all $i\in[n]$.
\end{algorithmic}
\caption{ Multi-agent A-MWU}
\end{algorithm}

\subsection{Parameters of A-MWU}
In practical applications, constant step is often used to obtain a quantative analysis of convergence rate. A necessary condition for the Riemannian gradient descent to avoid saddle points is that the step size $\alpha<\frac{1}{L}$, this condition ensures that the update rule is a loca diffeomorphism. In this section we give a discussion on the choice of the step size so that the Riemannian accelerated gradient descent can avoid saddle points provably. We derive a sufficient condition from the proof of Lemma \ref{DD} under which the RAGD can avoid saddle points. 

Recall that $\theta=\frac{s\gamma}{\gamma+s\mu}$ and $\zeta=\frac{(1-s)\gamma}{\bar{\gamma}}$, the following inequality gives a sufficient condition for the paramters:
\begin{equation}\label{ineq:step}
-b_d(c_d+1)+(c_d+1)^2>-b^2_dc_d
\end{equation}
where
\[
b_d=\left(\alpha(1-\theta)+\frac{s\theta}{\bar{\gamma}}\right)\lambda_d-\zeta(1-\theta)-1
\]
and 
\[
c_d=\zeta(1-\theta)(1-\alpha\lambda_d).
\]
Since $\lambda_d<0$ by assumption, then $b_d<0$. Thus a sufficient condition so that \ref{ineq:step} holds is $c_d>0$ which is equivalent to $\theta<1$ and then $s\gamma<\gamma+s\mu$. By the computation of $s$ and $\gamma$:
\[
s=\frac{\sqrt{\beta^2+4(1+\beta)\mu\alpha}-\beta}{2}
\] 
and
\[
\gamma=\frac{\sqrt{\beta^2+4(1+\beta)\mu\alpha}-\beta}{\sqrt{\beta^2+4(1+\beta)\mu\alpha}+\beta}.
\]

$s\gamma<\gamma+s\mu$ is computed as
\[
\frac{(\sqrt{\beta^2+4(1+\beta)\mu\alpha}-\beta)^2}{2(\sqrt{\beta^2+4(1+\beta)\mu\alpha}+\beta)}<\frac{\sqrt{\beta^2+4(1+\beta)\mu\alpha}-\beta}{\sqrt{\beta^2+4(1+\beta)\mu\alpha}+\beta}+\frac{\sqrt{\beta^2+4(1+\beta)\mu\alpha}-\beta}{2}\mu.
\]
After simplification, the above inequality is equivalent to 
\[
\beta^2+4(1+\beta)\mu\alpha<\left(\frac{2+(\mu+1)\beta}{1-\mu}\right)^2.
\]
$\mu$ can be taken to be small, i.e. smaller than 1 and the actual convexity parameter $\mu^*$, and then 
\[
\beta^2+4(1+\beta)\mu\alpha<(2+(\mu+1)\beta)^2
\]
sufficies. Furthermore, 
\[
\beta^2+4(1+\beta)\mu\alpha<4+4(\mu+1)\beta+(\mu+1)^2\beta^2
\]
and we note that this inequality can be implied by 
\[
4(1+\beta)\mu\alpha<4+4(\mu+1)\beta.
\]
Thus the step size satisfies
\[
\alpha<\frac{1+(\mu+1)\beta}{(1+\beta)\mu}.
\]
Therefore we have obtained a sufficient condition ensuring the saddle point avoidance that $\alpha$ should satisfy for any $\beta>0$. So in practice we can choose the constant step $\alpha<\min\{\frac{1}{L},\frac{1+(\mu+1)\beta}{(1+\beta)\mu}\}$, provided a chosen pair $\beta>0$ and $\mu<\min\{1,\mu^*\}$ where $\mu^*$ is the actual convexity parameter.

\subsection{Derivation of Accelerated Multiplicative Weights Update}
\subsubsection{MWU as Manifold Gradient Descent}
\paragraph{Orthogonality w.r.t. Shahshahani metric.}Since the positive orthant is a Riemannian manifold with the Shahshahani metric, the positive simplex, $\Delta^{d-1}_+=\{\vec{x}\in\mathbb{R}^d_+:\sum_jx_j=1\}$, has a natural submanifold structure by restricting the metric $g(\vec{x})$ on $\abs{\vec{x}}=1$, and thus $g_{ii}(\vec{x})=\frac{1}{x_i}$. The tangent space of $\Delta^{d-1}_+$ at $\vec{x}$ is denoted by $T_{\vec{x}}\Delta^{d-1}_+$. Note that $T_{\vec{x}}\Delta^{d-1}_+$ consists of all the vectors $\vec{v}$ satisfying $\sum_jv_j=0$, so we can identify all the tangent spaces on $\Delta^{d-1}_+$ with the hyperplane passing through $0$ and parallel to $\Delta^{d-1}_+$, i.e., 
\[
T_{\vec{x}}\Delta^{d-1}_+=\left\{\vec{v}\in\mathbb{R}^d:\sum_{j=0}^dv_j=0\right\}.
\]

\begin{figure}[H]
\centering
\includegraphics[width=0.7\textwidth]{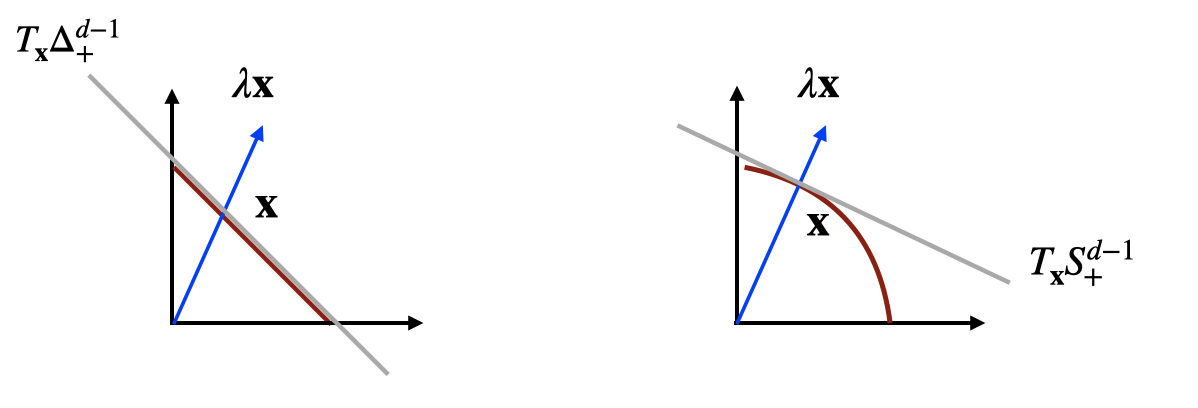}
\caption{figure}
\label{1}
\end{figure}

The most important feature of the Shahshahani manifold $\mathbb{R}^d_+$ is the way in which orthogonality is defined, especially the orthogonality with respect to $T_{\vec{x}}\Delta^{d-1}_+$ or $\Delta^{d-1}_+$. We denote $\langle\cdot,\cdot\rangle_{\vec{x}}$ the space-dependent inner product defined by the Shahshahani metric $g(\vec{x})$, and the following fact plays an essential role in comparing MWU to the manifold gradient descent on a standard sphere:
\begin{center}
\begin{tcolorbox}[enhanced,width=4.5in,center upper,
    fontupper=\bfseries,drop fuzzy shadow southwest,
    boxrule=0.4pt,sharp corners]
\emph{For all $\vec{u}\in T_{\vec{x}}\Delta^{d-1}_+$ and any $\lambda\ne 0$, it holds that $\langle\vec{u},\lambda\vec{x}\rangle_{\vec{x}}=0$.}
\end{tcolorbox}
\end{center}
%Now suppose $\vec{u}\in T_{\vec{x}}\Delta^{d-1}_+$, or equivalently $\vec{u}$ is parallel to the plane where the simplex lies. We are interested in the set of vectors $\vec{v}$ such that then the inner product between $\vec{u}$ and $\vec{v}$ is 
The verification is straightforward: $\langle\vec{u},\lambda\vec{x}\rangle_{\vec{x}}=\lambda\abs{\vec{x}}\sum_ju_j=0$ since $\sum_ju_j=0$. Geometrically, this statement means that the straight line passing through $0$ and $\vec{x}$ is orthogonal to the tangent space of $\Delta^{d-1}_+$ at $\vec{x}$, with respect to the Shahshahani metric on $\mathbb{R}^d_+$. This is very similar to the case in the Euclidean space where the line passing through the origin and a point $\vec{x}$ on a sphere is orthogonal to the tangent space at $\vec{x}$. See the illustration of Figure \ref{1} for the case of $\mathbb{R}^2$.

\begin{figure}[H]
\centering
\includegraphics[width=0.38\textwidth]{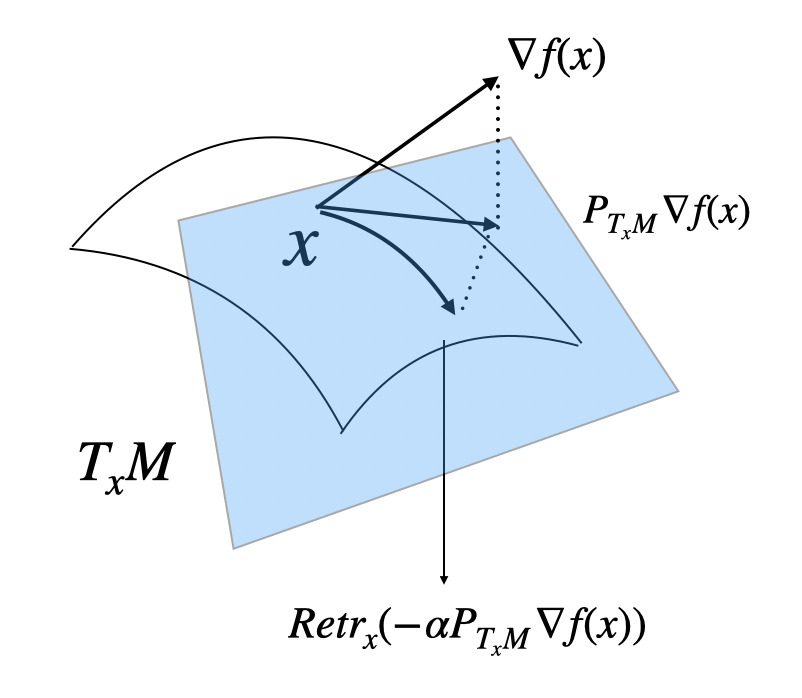}
\hspace{0.5in}
\includegraphics[width=0.38\textwidth]{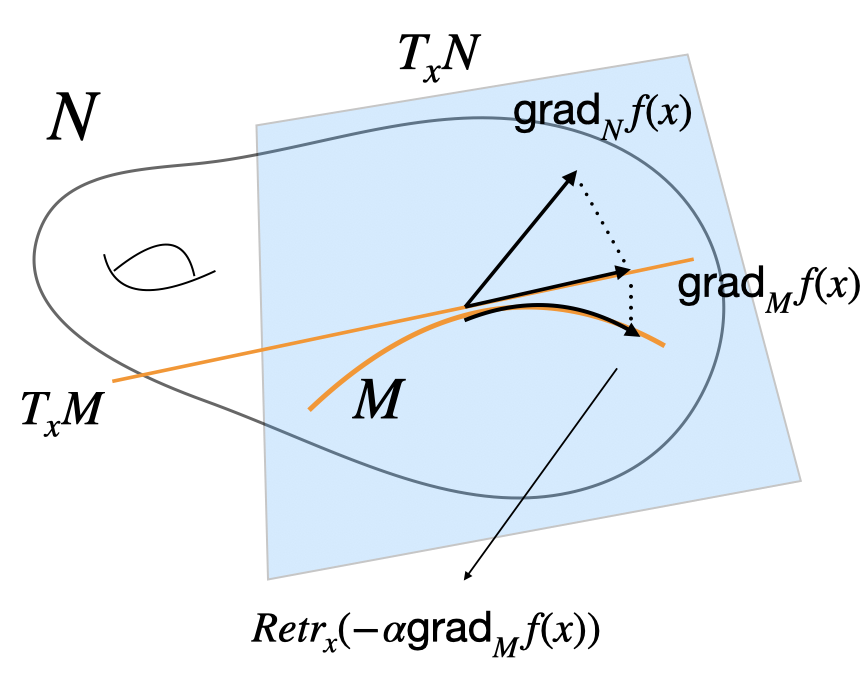}
\caption{RGD on $M$ embedded in $\mathbb{R}^d$ (left), and RGD on $M$ (the yellow curve) embedded in a manifold $N$ (right).}
\label{submnfldGD}
\end{figure}

 We next give a simply intuitive explanation that MWU is indeed a "Spherical gradient descent". By the manifold gradient descent for a function defined on $M\subset\mathbb{R}^d$ we mean the algorithm
\begin{equation}\label{RGD}
\vec{x}_{t+1}=\Retr_{\vec{x}_t}(-\alpha \mathcal{P}_{T_{\vec{x}}}\nabla f(\vec{x}_t))
\end{equation}
where $\mathcal{P}_{T_{\vec{x}}}\nabla f(\vec{x}_t)$ is the orthogonal projection of $\nabla f(\vec{x}_t)$ onto the tangent space $T_{\vec{x}_t}M$ with respect to the Euclidean metric on the ambient space $\mathbb{R}^d=T_{\vec{x}_t}\mathbb{R}^d$. To understand the intuition of MWU from a spherical viewpoint, we firstly need to generalize the manifold gradient descent (\ref{RGD}) to the case when the ambient space is a general Riemannian manifold $N$ instead of the Euclidean space $\mathbb{R}^d$. Denote $\grad_{N} f$ the Riemannian gradient of $f$ on $N$, then the generalized algorithm is written as
\[
\vec{x}_{t+1}=\Retr_{\vec{x}_t}\left(-\alpha\mathcal{P}_{T_{\vec{x}_t}M}\grad_Nf(\vec{x}_t)\right),
\]
where the orthogonal projection of $\grad_Nf(\vec{x})$ onto the tangent space $T_{\vec{x}}M$ is based on the inner product $\langle\cdot,\cdot\rangle_{\vec{x}}$ on $T_{\vec{x}}N$, i.e., the Riemannian metric on $N$. Figure \ref{submnfldGD} shows the intuition from gradient descent on $M\subset\mathbb{R}^d$ to that on $M\subset N$.

For the case of simplex, let $N=\mathbb{R}^d_+$ and $M=\Delta^{d-1}_+$, and the orthogonal projection in the tangent space is with respect to the Shahshahani metric. Then for small $\alpha>0$, the orthogonal projection of $-\alpha\grad_N f(\vec{x})$ onto $T_{\vec{x}}M$ is the vector obtained from the difference between the point $V(\vec{x})$ that is the normalization of $\vec{x}-\alpha \grad f(\vec{x})$ onto the simplex and the initial point $\vec{x}$, i.e.
\begin{align}
-\alpha\grad_Mf(\vec{x})&=-\alpha\mathcal{P}_{T_{\vec{x}}M}\grad_N f(x)
\\
&=\mathcal{P}_{T_{\vec{x}}M}\left(-\alpha\grad_N f(\vec{x})\right)
\\
&=V(\vec{x})-\vec{x}
\end{align}
where 
\[
V(\vec{x})=\left(\frac{x_1-\alpha x_1\frac{\partial f}{\partial x_1}}{1-\alpha\sum_jx_j\frac{\partial f}{\partial x_j}},...,\frac{x_d-\alpha x_d\frac{\partial f}{\partial x_d}}{1-\alpha\sum_jx_j\frac{\partial f}{\partial x_j}}\right),
\]
that is exactly the multiplicative weights update.
Thus by (\ref{RGD}), the multiplicative weights update is equivalent to the Riemannian gradient descent in Shahshahani manifold in that
\begin{align}
\Retr_{\vec{x}_t}(-\alpha\grad_Mf(\vec{x}_t))&=\vec{x}-\alpha\grad_Mf(\vec{x}_t)=\vec{x}_t+V(\vec{x}_t)-x_t=V(\vec{x}_t).
\end{align}

\paragraph{Derivation of A-MWU.} The exponential map on the positive simplex with Shahshahani metric is a well known result. Given a point $\vec{x}\in\Delta^{d-1}_+$ and a vector $\vec{v}\in T_{\vec{x}}\Delta^{d-1}_+$, the exponential map is
\[
\Exp_{\vec{x}}(\vec{u})=\left(\frac{x_1e^{u_1}}{\sum_jx_je^{u_j}},...,\frac{x_de^{u_d}}{\sum_jx_je^{u_j}}\right).
\]
Then the logarithmic map can be computed explicitly, i.e., given a base point $\vec{x}$ and any point $\vec{y}$ on $\Delta^{d-1}_+$, we need to find the expression of the vector $\vec{u}\in T_{\vec{x}}\Delta^{d-1}_+$ such that 
\[
\Exp_{\vec{x}}(\vec{u})=\vec{y}.
\]
 Let 
 \[
 S=\sum_jx_je^{u_j},
 \]
  then the vector equation 
  \[
  \Exp_{\vec{x}}(\vec{u})=\vec{y}
  \]
   gives a system of equations, for $i\in[d]$, there is
\[
\frac{x_ie^{u_i}}{S}=y_i.
\]
Multiplying them up, we have
\[
\prod_{i=1}^d\frac{x_ie^{u_i}}{S}=\frac{1}{S^d}\left(\prod_{i=1}^d x_i\right)\left(\prod_{i=1}^d e^{u_i}\right)=\prod_{i=1}^dy_i.
\]
Since $\vec{u}=(u_1,...,u_d)\in T_{\vec{x}}\Delta^{d-1}_+$, we have $\sum_ju_j=0$ and thus $\prod_{i=1}^de^{u_j}=e^{\sum_ju_j}=1$. Then the normalization constant $S$ can be computed as
$
S=\left(\prod_{i=1}^d\frac{x_i}{y_i}\right)^{1/d}.
$
Moreover, the equations $\frac{x_ie^{u_i}}{S}=y_i$ gives the solution of $u_i$ by
$
u_i=\ln\left(S\frac{y_i}{x_i}\right)
$, and then the logarithmic map is 
\[
\Log_{\vec{x}}(\vec{y})=\left(\ln\left(S\frac{y_1}{x_1}\right),...,\ln\left(S\frac{y_d}{x_d}\right)\right).
\]

\subsection{Missing Proofs}

\subsubsection{Proof of Theorem \ref{generalSMT}.}

Before going through the details of the proof, we first give a sketched idea and review some qualitative theory of Ordinary Differential Equations.

The main technique of proving the above result is the \textbf{Discrete Lyapunov-Perron Method} for non-autonomous dynamical systems. Historically the Lyapunov-Perron method is often used to prove the existence of stable manifold for continuous time dynamical system defined by ordinary differential equations of the form
\[
\frac{d\vec{x}}{dt}=A(t)\vec{x}+R(t,\vec{x}),
\] 
where $A(t)$ is a (time-dependent) matrix. The idea is to write the general solution in the integral form and deduce that the solution $u(t,\vec{x}_0)$ with initial condition $\vec{x}_0$ converging to an unstable fixed point must be the fixed point of the integral operator $T$ defined on the space of continuous functions as
\begin{align}
Tu(t,\vec{x}_0)&=U(t)\vec{x}_0+\int_0^tU(t-s)R(s,u(s,\vec{x}_0))ds-\int_t^{\infty}V(t-s)R(s,u(s,\vec{x}_0))ds.
\end{align}
The discrete Lyapunov-Perron method for the dynamical system (\ref{eq:SMT1}) results in an operator $T$ defined on the space of sequences (in contrast to continuous time functions), which transforms a sequence $\{\vec{x}_t\}_{t\in\mathbb{N}}$ into a new one $\{(T\vec{x})_t\}_{t\in\mathbb{N}}$ in the following way:
\begin{align}\label{eq:T}
(T\vec{x})_{t+1}=\left(B(t,0)\vec{x}_0^++\sum_{i=0}^tB(t,i+1)\vect{\eta}^+(i,\vec{x}_i)\right)
\oplus\left(-\sum_{i=0}^{\infty}C(t+1+i,t+1)^{-1}\vect{\eta}^-(t+1+i,\vec{x}_{t+1+i})\right),
\end{align}
where ``$\oplus$" refers to the stable-unstable decomposition $E_s\oplus E_u$.
Another ingredient show the existence of the local stable manifold is the \textbf{Banach Fixed Point Theorem}, \cite{AL}, saying that a contraction map $T:\mathcal{V}\rightarrow\mathcal{V}$ on a complete metric space $\mathcal{V}$ has unique fixed point. By a contraction map, we mean $d(Tx,Ty)\le Kd(x,y)$ for some $K<1$. This ensures the existence and uniqueness of local stable manifold of the dynamical system (\ref{eq:SMT1}).

\begin{proof}

Denote $A(m,n)$ the matrix as follows
\[
A(m,n)=\left[
\begin{array}{ccc}
\prod_{i=n}^m\mathcal{L}_1(i)&&
\\
&\ddots&
\\
&&\prod_{i=n}^m\mathcal{L}_d(i)
\end{array}
\right]
\]
Since the stable and unstable subspaces are only determined by the signs of eigenvalue of the Hessian $\nabla^2f(\vec{y})$, the matrix $A(m,n)$ can be decomposed to stable-unstable blocks, i.e.
\[
A(m,n)=\left[
\begin{array}{cc}
B(m,n)&
\\
&C(m,n)
\end{array}
\right]
\]
such that there exit positive numbers $K_1,K_2<1$ satisfying 
\begin{align}
\norm{B(m,n)}_{sp}&\le K_1^{m-n+1}
\\
\norm{C(m,n)^{-1}}_{sp}&\le K_2^{m-n+1}.
\end{align}

Starting from the intial condition $\vec{x}_0$, the dynamical system can be represented by
\[
\vec{x}_{t+1}=A(t,0)\vec{x}_0+\sum_{i=0}^tA(t,i+1)\vect{\eta}(i,\vec{x}_i)
\]
which can be decomposed into stable-unstable parts 
\[
\vec{x}_{t+1}^+=B(t,0)\vec{x}_0^++\sum_{i=0}^tB(t,i+1)\vect{\eta}^+(i,\vec{x}_i)
\]
\[
\vec{x}_{t+1}^-=C(t,0)\vec{x}_0^-+\sum_{i=0}^t C(t,i+1)\vect{\eta}^-(i,\vec{x}_i).
\]

If $\vec{x}_{t+1}\rightarrow 0$ as $t\rightarrow\infty$, then $\vec{x}^-_{t+1}\rightarrow 0$ as $k\rightarrow\infty$. So the unstable component $\vec{x}_0^-$ of the initial point $\vec{x}_0$ satisfies 
\[
\vec{x}_0^-=-\sum_{i=1}^{\infty}C(i-1,0)^{-1}\vect{\eta}^-(i-1,\vec{x}_{i-1}),
\]
where $C(m,n)^{-1}$ denotes the inverse of $C(m,n)$. Therefore we can write the updated term $\vec{x}_{t+1}$ as

\begin{align}
\vec{x}_{t+1}&=\vec{x}_{t+1}^+\oplus \vec{x}_{t+1}^-
\\
&=\left(B(t,0)\vec{x}_0^++\sum_{i=0}^tB(t,i+1)\vect{\eta}^+(i,\vec{x}_i)\right)\oplus\left(C(t,0)\vec{x}_0^-+\sum_{i=0}^tC(t,i+1)\vect{\eta}^-(i,\vec{x}_i)\right)
\\
&=\left(B(t,0)\vec{x}_0^+\right)\oplus\left(-\sum_{i=0}^tC(t,i+1)\vect{\eta}^-(i,\vec{x}_i)-\sum_{i=0}^{\infty}C(t+1+i,t+1)^{-1}\vect{\eta}^-(t+1+i,\vec{x}_{k+1+i})\right)
\\
&+\left(\sum_{i=0}^tB(t,i+1)\vect{\eta}^+(i,\vec{x}_i)\right)\oplus\left(\sum_{i=0}^tC(t,i+1)\vect{\eta}^-(i,\vec{x}_i)\right)
\\
&=\left(B(t,0)\vec{x}_0^++\sum_{i=0}^tB(t,i+1)\vect{\eta}^+(i,\vec{x}_i)\right)\oplus\left(-\sum_{i=0}^{\infty}C(t+1+i,t+1)^{-1}\vect{\eta}^-(t+1+i,\vec{x}_{k+1+i})\right)
\end{align}

Define the integral operator acting on the space of bounded sequences as 
\[
(T\vec{x})_{t+1}=\left(B(t,0)\vec{x}_0^++\sum_{i=0}^tB(t,i+1)\vect{\eta}^+(i,\vec{x}_i)\right)\oplus\left(-\sum_{i=0}^{\infty}C(t+1+i,t+1)^{-1}\vect{\eta}^-(t+1+i,\vec{x}_{k+1+i})\right)
\]

To show the existence of local stable manfiold, the main step is to show that $T$ acting on the space of sequences is a contraction map if all the sequences considered are in a small Euclidean ball around $\vec{0}$. Thus we investigate the norm difference between the image of two sequences $\vec{u}=\{\vec{u}\}_{n\in\mathbb{N}}$ and $\vec{v}=\{\vec{v}_n\}_{n\in\mathbb{N}}$.

\begin{align}
(T\vec{u}-T\vec{v})_{t+1}&=(T\vec{u})_{t+1}-(T\vec{v})_{t+1}
\\
&=\left(B(t,0)(\vec{u}_0^+-\vec{v}_0^+)+\sum_{i=0}^kB(t,i+1)(\vect{\eta}^+(i,\vec{u}_i)-\vect{\eta}^+(i,\vec{i}))\right)
\\
&\oplus\left(-\sum_{i=0}^{\infty}C(t+1+i,t+1)^{-1}(\vect{\eta}^{-}(t+1+i,\vec{u}_{t+1+i})-\vect{\eta}^-(t+1+i,\vec{v}_{t+1+i}))\right)
%\\
%&=\left[
%\begin{array}{c}
%B(t,0)(\vec{u}_0^+-\vec{v}_0^+)+\sum_{i=0}^kB(t,i+1)(\vect{\eta}^+(i,\vec{u}_i)-\vect{\eta}^+(i,\vec{i}))
%\\
%-\sum_{i=0}^{\infty}C(t+1+i,t+1)^{-1}(\vect{\eta}^{-}(t+1+i,\vec{u}_{t+1+i})-\vect{\eta}^-(t+1+i,\vec{v}_{t+1+i}))
%\end{array}
%\right]
\end{align}

and then

\begin{align}
\norm{(T\vec{u}-T\vec{v})_{t+1}}&\le \norm{B(t,0)}_{sp}\norm{\vec{u}_0^+-\vec{v}_0^+}+\sum_{i=0}^t\norm{B(t,i+1)}_{sp}\norm{\vect{\eta}^+(i,\vec{u}_i)-\vect{\eta}^+(i,\vec{v}_i)}
\\
&+\sum_{i=0}^{\infty}\norm{C(t+1+i,t+1)^{-1}}_{sp}\norm{\vect{\eta}^-(t+1+i,\vec{u}_{t+1+i})-\vect{\eta}^-(t+1+i,\vec{v}_{t+1+i})}
\\
&\le\norm{B(t,0)}_{sp}\norm{\vec{u}_0^+-\vec{v}_0^+}+\sum_{i=0}^t\norm{B(t,i+1)}_{sp}\alpha_i\epsilon\norm{\vec{u}_i-\vec{v}_i}
\\
&+\sum_{i=0}^{\infty}\norm{C(k+1+i,k+1)^{-1}}_{sp}\alpha_{t+1+i}\epsilon\norm{\vec{u}_{t+1+i}-\vec{v}_{t+1+i}}
\\
&\le \norm{B(t,0)}d(\vec{u},\vec{v})+\sum_{i=0}^t\norm{B(t,i+1)}\alpha_i\epsilon d(\vec{u},\vec{v})
\\
&+\sum_{i=0}^{\infty}\norm{C(t+1+i,t+1)^{-1}}_{sp}\alpha_{t+1+i}\epsilon d(\vec{u},\vec{v})
\\
&=\left(\norm{B(t,0)}_{sp}+\epsilon\left(\sum_{i=0}^t\alpha_i\norm{B(t,i+1)}_{sp}\right)+\epsilon\left(\sum_{i=0}^{\infty}\alpha_{t+1+i}\norm{C(t+1+i,t+1)^{-1}}_{sp}\right)\right)d(\vec{u},\vec{v}).
\end{align}
By the following lemma \ref{lemma:contraction}, we conclude that $T$ is a contraction map. By Banach fixed point theorem, there exists a unique sequence $\vec{x}=\{\vec{x}_n\}_{n\in\mathbb{N}}$, such that
\[
T\vec{x}=\vec{x}.
\]
For this sequence, the initial point $\vec{x}_0$ should satisfy 
\begin{equation}\label{eq:x0}
(\vec{x}_0^+,\vec{x}_0^-)=(\vec{x}_0^+,-\sum_{i=0}^{\infty}C(t+1+i,t+1)\vect{\eta}^-(t+1+i,\vec{x}_{t+1+i})).
\end{equation}
Since the sequence is generated by the dynamical system with a specific initial point $\vec{x}_0$, each term $\vec{x}_t$ is determined completely by the initial point. Thus we can consider each term $\vec{x}_t$ as a function of the initial point, i.e. 
\[
\vec{x}_t=\vec{x}_t(\vec{x}_0)=\vec{x}_t(\vec{x}_0^+,\vec{x}_0^-).
\]
Pluggin the above form to the equality \ref{eq:x0}, we have
\[
(\vec{x}_0^+,\vec{x}_0^-)=(\vec{x}_0^+,-\sum_{i=0}^{\infty}C(t+1+i,t+1)\vect{\eta}^-(t+1+i,\vec{x}_{t+1+i}(\vec{x}_0^+,\vec{x}_0^-)))
\]
and the unstable component gives the following equation of $\vec{x}_0^+$ and $\vec{x}_0^-$,
\begin{equation}\label{eq:x00}
\vec{x}_0^-=-\sum_{i=0}^{\infty}C(t+1+i,t+1)\vect{\eta}^-(t+1+i,\vec{x}_{t+1+i}(\vec{x}_0^+,\vec{x}_0^-)).
\end{equation}
With a specific dynamical system, the right hand side is completely determined by $\vec{x}_0^+$ and $\vec{x}_0^-)$, so we call the infinite sum on the right hand side a function $\Phi(\vec{x}_0^+,\vec{x}_0^-)$, i.e.
\[
\Phi(\vec{x}_0^+,\vec{x}_0^-)=-\sum_{i=0}^{\infty}C(t+1+i,t+1)\vect{\eta}^-(t+1+i,\vec{x}_{t+1+i}(\vec{x}_0^+,\vec{x}_0^-))
\]
and equation \ref{eq:x00} is written as
\[
\vec{x}_0^-=\Phi(\vec{x}_0^+,\vec{x}_0^-).
\]
This equation defines an implicit function $\vec{x}_0^-=\varphi(\vec{x}_0^+)$ by the uniqueness of Banach fixed point theorem. The stable manifold is nothing but the graph of $\varphi$.
\end{proof}

%\[
%(T\vec{x})_{t+1}=\left[
%\begin{array}{c}
%B(t,0)\vec{x}_0^++\sum_{i=0}^tB(t,i+1)\vect{\eta}^+(i,\vec{x}_i)
%\\
%-\sum_{i=0}^{\infty}C(t+1+i,t+1)^{-1}\vect{\eta}^-(t+1+i,\vec{x}_{k+1+i})
%\end{array}
%\right]
%\]

\begin{lemma}\label{lemma:contraction}
$T$ is a contraction map.
\end{lemma}

\begin{proof}
Let $\vec{u}=\{\vec{u}\}_{n\in \mathbb{N}}$ and $\vec{u}=\{\vec{v}_n\}_{n\in\mathbb{N}}$ be two sequences. Based on the assumptions of $\vect{\eta}$, we have
\begin{align}
&\norm{(T\vec{u}-T\vec{v})_{t+1}}
\\
&\le\left(\norm{B(t,0)}_{sp}+\epsilon\left(\sum_{i=0}^t\norm{B(t,i+1)}_{sp}\right)+\epsilon\left(\sum_{i=0}^{\infty}\norm{C(t+1+i,t+1)^{-1}}_{sp}\right)\right)d(\vec{u},\vec{v})
\end{align}
 By choosing $\epsilon>0$ properly, one can make $T$ a contraction map. To be specific, we compute
 \begin{align}
 \sum_{i=0}^t\norm{B(t,i+1)}_{sp}\le\sum_{i=0}^tK_1^{t-i}
 =\frac{1-K_1^{t+1}}{1-K_1}\le\frac{1}{1-K_1}
 \end{align}
 \begin{align}
 \sum_{i=0}^{\infty}\norm{C(t+1+i,t+1)^{-1}}_{sp}\le\sum_{i=0}^{\infty}K_2^{i+1}&=\lim_{N\rightarrow\infty}\sum_{i=0}^NK_2^{i+1}
 \\
 &=\lim_{N\rightarrow\infty}\frac{K_2(1-K_2^{N+1})}{1-K_2}
 \\
 &=\frac{K_2}{1-K_2}
 \end{align}
 \[
 \norm{B(t,0)}_{sp}\le K_1^t\le K_1<1.
 \]
 So we need to choose $\epsilon$ small such that
 \[
 \norm{B(t,0)}_{sp}+\epsilon\frac{1}{1-K_1}+\epsilon\frac{K_2}{1-K_2}\le1.
 \]
 For example, take $\epsilon=\frac{1}{2}\frac{(1-K_1)^2(1-K_2)}{1-K_1K_2}$ and then the operator $T$ is a contraction map.
 \end{proof}
 
\subsection{Proof of Proposition \ref{globalmeasurezero}}

\begin{proof}
The proof follows \cite{PPW19a}. Throughout this proof, we define $\tilde{\psi}(m,n,\vec{x})=\psi(m,...,\psi(n+1,\psi(n,\vec{x}))...)$ for $m>n$.
 For each $\vec{x}^*\in \mathcal{A}^*$, there is an associated open neighborhood $U_{\vec{x}^*}$ %promised by Theorem \ref{Center-Stable Manifold Theorem}, \ref{SMT:local} or \ref{generalSMT}, 
 depending on the dynamical system we consider. $\bigcup_{\vec{x}^*\in \mathcal{A}^*}U_{\vec{x}^*}$ forms an open cover, and since manifold $M$ is second-countable, we can find a countable subcover, so that 
\[\bigcup_{\vec{x}^*\in \mathcal{A}^*}U_{\vec{x}^*}=\bigcup_{i=1}^{\infty}U_{\vec{x}^*_i}.
\]
By the definition of global stable set, we have
\[
W^s(\mathcal{A}^*)=\{\vec{x}_0:\lim_{k\rightarrow\infty}\tilde{\psi}(k,0,\vec{x}_0)\in\mathcal{A}^*\}.
\]
Fix a point $\vec{x}_0\in W^s(\mathcal{A}^*)$. Since 
\[
\tilde{\psi}(k,0,\vec{x}_0)\rightarrow\vec{x}^*\in\mathcal{A}^*,
\]
 there exists some nonnegative integer $T$ and all $t\ge T$, such that
\[
\tilde{\psi}(t,0,\vec{x}_0)\in\bigcup_{\vec{x}^*\in\mathcal{A}^*}U_{\vec{x}^*}=\bigcup_{i=1}^{\infty}U_{\vec{x}^*_i}.
\]
So $\tilde{\psi}(t,0,\vec{x}_0)\in U_{\vec{x}^*_i}$ for some $\vec{x}^*_i\in\mathcal{A}^*$ and all $t\ge T$. This is equivalent to
\[
\tilde{\psi}(T+k,T,\tilde{\psi}(T,0,\vec{x}_0))\in U_{\vec{x}^*_i}
\]
for all $k\ge 0$, and this implies that
\[
\tilde{\psi}(T,0,\vec{x}_0)\in\tilde{\psi}^{-1}(T+k,T,U_{\vec{x}^*_i})
\]
for all $k\ge 0$. And then we have
\[
\tilde{\psi}(T,0,\vec{x}_0)\in\bigcap_{k=0}^{\infty}\tilde{\psi}^{-1}(T+k,T,U_{\vec{x}^*_i}).
\]
Denote $S_{i,T}:=\bigcap_{k=0}^{\infty}\tilde{\psi}^{-1}(T+k,T,U_{\vec{x}^*_i})$
and the above relation is equivalent to
$\vec{x}_0\in \tilde{\psi}^{-1}(T,0,S_{i,T})$.
Take the union for all nonnegative integers $T$, we have
\[
\vec{x}_0\in\bigcup_{T=0}^{\infty}\tilde{\psi}^{-1}(T,0,S_{i,T}).
\]
And union for all $i$ we obtain that
\[
\vec{x}_0\in \bigcup_{i=1}^{\infty}\bigcup_{T=0}^{\infty}\tilde{\psi}^{-1}(T,0,S_{i,T})
\]
implying that
\[
W^s(\mathcal{A}^*)\subset\bigcup_{i=1}^{\infty}\bigcup_{T=0}^{\infty}\tilde{\psi}^{-1}(T,0,S_{i,T}).
\]
Since $S_{i,T}\subset W_n(\vec{x}^*)$, and $W_n(\vec{x}^*)$ has codimension at least 1. This implies that $S_{i,T}$ has measure 0 with respect to the volume measure from the Riemannian metric on $M$. Since the image of set of measure zero under diffeomorphism is of measure zero, and countable union of zero measure sets is still measure zero, we obtain that $W^s(\mathcal{A}^*)$ is of measure zero.
\end{proof}

\subsection{Proof of Theorem \ref{SMT:RAGD2}}

\begin{proposition}\label{exp:taylor}
The differential of $\Exp_{\vec{x}}(\vec{v})$ at $\vec{v}=0\in T_{\vec{x}}M$ is the identity map, i.e.
\[
D\Exp_{\vec{x}}(\vec{0})=Id.
\]
Moreover, the exponential map has the expansion 
\[
\Exp_{\vec{x}}(\vec{v})=\vec{x}+\vec{v}-\frac{1}{2}\sum_{i,j}\Gamma_{ij}^kv_iv_j+O(\norm{\vec{v}}^3)
\]
\end{proposition}

\begin{proof}
\begin{align}
D\Exp_{\vec{x}}(\vec{0})[\vec{v}]&=\frac{d}{dt}(\Exp_{\vec{x}}(t\vec{v}))\Big|_{t=0}
\\
&=\frac{d}{dt}(\gamma(1,\vec{x},t\vec{v}))\Big|_{t=0}
\\
&=\frac{d}{dt}(\gamma(t,\vec{x},\vec{v}))\Big|_{t=0}
\\
&=\vec{v}.
\end{align}

Since the geodesic equation in local coordinate system is written as
\[
\frac{d^2x_k}{dt^2}+\sum_{i,j}\Gamma^k_{ij}\frac{dx_i}{dt}\frac{dx_j}{dt}=0,
\]
therefore,
\[
\Exp_{\vec{x}}(\vec{v})=\vec{x}+\vec{v}-\frac{1}{2}\sum_{i,j}\Gamma_{ij}^kv_iv_j+O(\norm{\vec{v}}^3).
\]
\end{proof}

\begin{remark}
From the first order approximation of the exponential map, we have the approximation of the logarithmic map on manifold,
\[
\Log_{\vec{x}}(\vec{z})=\vec{z}-\vec{x}+O(d^2(\vec{x},\vec{z}))
\]
where $d(\vec{x},\vec{z})$ is the geodesic distance between $\vec{x}$ and $\vec{z}$ on $M$.
\end{remark}

The following lemma plays an essential role in proving the saddle avoidance result.
\begin{lemma}[Eigenvalue]\label{lemma:eigenvalue}
The map $\psi$ given by RAGD is a local diffeomorphism on $M\times M$. Moreover, if $\vec{x}^*$ is a saddle point of $f:M\rightarrow\mathbb{R}$, then $(\vec{x}^*,\vec{x}^*)$ is an unstable fixed point of $\psi$.
\end{lemma}
\begin{proof}
First of all we give a sketch of the strategy of proving this lemma.
To show $(\vec{x}^*,\vec{x}^*)$ is an unstable fixed point of $\psi$, it suffices to show that the differential $D\psi(t,\vec{x}^*,\vec{x}^*)$ has an eigenvalue with magnitude greater than 1. Since the direct computation shows that the Jacobian $D\psi(t,\vec{x}^*,\vec{x}^*)$ at $(\vec{x}^*,\vec{x}^*)$ is of the form:
\[
D\psi(t,\vec{x}^*,\vec{x}^*)=
\left[
\begin{array}{cc}
D_{\vec{x}}F& D_{\vec{v}}F
\\
D_{\vec{x}}G & D_{\vec{v}}G
\end{array}
\right]_{\vec{x}^*}
\]
where the blocked matrices evaluated at $\vec{x}^*$ as follows
\begin{align}
D_{\vec{x}}F&=(1-\theta_t)(I-\alpha_t D_{\vec{y}}\grad f(\vec{y}))|_{\vec{y}=\vec{x}^*}
\\
D_{\vec{v}}F&=\theta_t(I-\alpha_t D_{\vec{y}}\grad f(\vec{y}))|_{\vec{y}=\vec{x}^*}
\\
D_{\vec{x}}G&=\left(1-\theta_t\right)\left((1-\zeta_t)I-\frac{s_t}{\bar{\gamma}_t}D_{\vec{y}}\grad f(\vec{y})\right)\Big|_{\vec{y}=\vec{x}^*}
\\
D_{\vec{v}}G&=((1-\zeta_t)\theta_t+\zeta_t)I-\frac{s_t\theta_t}{\bar{\gamma}_t}D_{\vec{y}}\grad f(\vec{y})|_{\vec{y}=\vec{x}^*}
\end{align}
and
$\theta_t=\frac{s_t\gamma_t}{\gamma_t+s_t\mu}$, $\zeta_t=\frac{(1-s_t)\gamma_t}{\bar{\gamma}_t}$. By choosing local coordinate properly, e.g. the normal coordinate, we can make the differential of the Riemannian gradient $D_{\vec{y}}\grad f(\vec{y})$ equal to the Euclidean Hessian $\nabla^2f(\vec{y})$ which is similar to a diagonal matrix. Therefore the deteminant of diagonalized $D\psi(t,\vec{x}^*,\vec{x}^*)-xI$ can be computed explicitly and the existence of eigenvalue whose magnitude greater than 1 follows.
Now let's make the argument complete.

Recall that with convex parameter $\mu$, the Riemannian accelerated gradient descent algorithm is written as follows:
\[
\vec{y}_t=\Exp_{\vec{x}_t}\left(\frac{s_t\gamma_t}{\gamma_t+s_t\mu}\Log_{\vec{x}_t}(\vec{v}_t)\right)
\]

\[
\vec{x}_{t+1}=\Exp_{\vec{y}_t}(-\alpha_t\grad f(\vec{y}_t))
\]

\[
\vec{v}_{t+1}=\Exp_{\vec{y}_t}\left(\frac{(1-s_t)\gamma_t}{\bar{\gamma}_t}\Log_{\vec{y}_t}(\vec{v}_t)-\frac{s_t}{\bar{\gamma}_t}\grad f(\vec{y}_t)\right)
\]
We note that if the step size is taken to be 0, the above algorithm is an identity map on the manifold $M$. So the determinant of the differential of the above map goes to 1 as $\alpha_t$ goes to 0. Since the determinant is a continuous function with respect to its entries, by choosing small enough step size, the determinant of the differential is positive and this implies the above algorithm is a local diffeomorphism.

Locally, these expression has the expansion by Proposition \ref{exp:taylor},
\begin{align}
\vec{y}_t&=\Exp_{\vec{x}_t}\left(\frac{s_t\gamma_t}{\gamma_t+s_t\mu}\Log_{\vec{x}_t}(\vec{v}_t)\right)
\\
&=\vec{x}_t+\frac{s_t\gamma_t}{\gamma_t+s_t\mu}\Log_{\vec{x}_t}(\vec{v}_t)+O\left(\norm{\frac{s_t\gamma_t}{\gamma_t+s_t\mu}\Log_{\vec{x}_t}(\vec{v}_t)}^2\right)
\\
&=\vec{x}_t+\frac{s_t\gamma_t}{\gamma_t+s_t\mu}(\vec{v}_t-\vec{x}_t+O(d^2(\vec{x}_t,\vec{v}_t)))+O\left(\norm{\frac{s_t\gamma_t}{\gamma_t+s_t\mu}\Log_{\vec{x}_t}(\vec{v}_t)}^2\right)
\\
&=\left(1-\frac{s_t\gamma_t}{\gamma_t+s_t\mu}\right)\vec{x}_t+\frac{s_t\gamma_t}{\gamma_t+s_t\mu}\vec{v}_t+\frac{s_t\gamma_t}{\gamma_t+s_t\mu}O(d^2(\vec{x}_t,\vec{v}_t))+O\left(\norm{\frac{s_t\gamma_t}{\gamma_t+s_t\mu}\Log_{\vec{x}_t}(\vec{v}_t)}^2\right)
\\
&=\left(1-\frac{s_t\gamma_t}{\gamma_t+s_t\mu}\right)\vec{x}_t+\frac{s_t\gamma_t}{\gamma_t+s_t\mu}\vec{v}_t+O(d^2(\vec{x}_t,\vec{v}_t))
\end{align}

\begin{align}
\vec{x}_{t+1}&=\Exp_{\vec{y}_t}(-\alpha_t\grad f(\vec{y}_t))
\\
&=\vec{y}_t-\alpha_t\grad f(\vec{y}_t)+O(\norm{\alpha_t\grad f(\vec{y}_t)}^2)
\end{align}

\begin{align}
\vec{v}_{t+1}&=\Exp_{\vec{y}_t}\left(\frac{(1-s_t)\gamma_t}{\bar{\gamma}_t}\Log_{\vec{y}_t}(\vec{v}_t)-\frac{s_t}{\bar{\gamma}_t}\grad f(\vec{y}_t)\right)
\\
&=\vec{y}_t+\frac{(1-s_t)\gamma_t}{\bar{\gamma}_t}\Log_{\vec{y}_t}(\vec{v}_t)-\frac{s_t}{\bar{\gamma}_t}\grad f(\vec{y}_t)+O(\norm{V}^2)
\\
&=\vec{y}_t+\frac{(1-s_t)\gamma_t}{\bar{\gamma}_t}(\vec{v}_t-\vec{y}_t+O(d^2(\vec{y}_t,\vec{v}_t)))-\frac{s_t}{\bar{\gamma}_t}\grad f(\vec{y}_t)+O(\norm{V}^2)
\\
&=\left(1-\frac{(1-s_t)\gamma_t}{\bar{\gamma}_t}\right)\vec{y}_t+\frac{(1-s_t)\gamma_t}{\bar{\gamma}_t}\vec{v}_t-\frac{s_t}{\bar{\gamma}_t}\grad f(\vec{y}_t)+O(\norm{V}^2)
\end{align}
where
\[
V=\frac{(1-s_t)\gamma_t}{\bar{\gamma}_t}\Log_{\vec{y}_t}(\vec{v}_t)-\frac{s_t}{\bar{\gamma}_t}\grad f(\vec{y}_t).
\]
It is easy to see that as $t\rightarrow\infty$, the sequences of $\{\vec{x}_t\}_{t\in\mathbb{N}}$, $\{\vec{y}_t\}_{t\in\mathbb{N}}$ and $\{\vec{v}_t\}_{t\in\mathbb{N}}$ generated by the algorithm converge to the same point $\vec{x}^*$, the critical point of $f(\vec{x})$.

Let $\theta_t=\frac{s_t\gamma_t}{\gamma_t+s_t\mu}$ and $\zeta_t=\frac{(1-s_t)\gamma_t}{\bar{\gamma}_t}$, we denote
\begin{align}
F(\vec{y})&=\vec{y}-\alpha_t\grad f(\vec{y})+O(\norm{\alpha_t\grad f(\vec{y})}^2)
\\
G(\vec{y},\vec{v})&=\left(1-\zeta_t\right)\vec{y}_t+\zeta_t\vec{v}_t-\frac{s_t}{\bar{\gamma}_t}\grad f(\vec{y}_t)+O(\norm{V}^2)
\\
\vec{y}(\vec{x},\vec{v})&=(1-\theta_t)\vec{x}+\theta_t\vec{v}+O(d^2(\vec{x},\vec{v}))
\end{align}
and then the algorithm can be written as
\begin{align}
(\vec{x},\vec{v})\leftarrow\psi(t,\vec{x},\vec{v})=(F(\vec{y}),G(\vec{y},\vec{v}))
\end{align}

The differential of $\psi$:

\[
D\psi(t,\vec{x},\vec{v})=
\left[
\begin{array}{cc}
D_{\vec{x}}F&D_{\vec{v}}F
\\
D_{\vec{x}}G&D_{\vec{v}}G
\end{array}
\right]
\]

\begin{align}
D_{\vec{x}}F&=D_{\vec{y}}F\circ D_{\vec{x}}\vec{y}
\\
&=\left(I-\alpha_t D_{\vec{y}}\grad f(\vec{y})+O(\norm{\alpha_t\grad f(\vec{y})})\right)\circ\left((1-\theta_t)I+O(d(\vec{x},\vec{v}))\right)
\\
&=(1-\theta_t)(I-\alpha_t D_{\vec{y}}\grad f(\vec{y}))+(1-\theta_t)O(\norm{\alpha_t\grad f(\vec{y})})
\\
&+(I-\alpha_t D_{\vec{y}}\grad f(\vec{y}))\circ O(d(\vec{x},\vec{v}))+O(\norm{\alpha_t\grad f(\vec{y})})\circ O(d(\vec{x},\vec{v}))
\end{align}

\begin{align}
D_{\vec{v}}F&=D_{\vec{y}}F\circ D_{\vec{v}}\vec{y}
\\
&=\left(I-\alpha_t D_{\vec{y}}\grad f(\vec{y})+O(\norm{\alpha_t\grad f(\vec{y})})\right)\circ(\theta_t I+O(d(\vec{x},\vec{v})))
\\
&=\theta_t(I-\alpha_t D_{\vec{y}}\grad f(\vec{y}))+\theta_t O(\norm{\alpha_t\grad f(\vec{y})})
\\
&+(I-\alpha_t D_{\vec{y}}\grad f(\vec{y}))\circ O(d(\vec{x},\vec{y}))+O(\norm{\alpha_t\grad f(\vec{y})})\circ O(d(\vec{x},\vec{v}))
\end{align}

\begin{align}
D_{\vec{x}}G&=D_{\vec{y}}G\circ D_{\vec{x}}\vec{y}
\\
&=\left((1-\zeta_t)I-\frac{s_t}{\bar{\gamma}_t}D_{\vec{y}}\grad f(\vec{y})+O(\norm{V})\right)\circ((1-\theta_t)I+O(d(\vec{x},\vec{y})))
\\
&=(1-\theta_t)\left((1-\zeta_t)I-\frac{s_t}{\bar{\gamma}_t}D_{\vec{y}}\grad f(\vec{y})\right)+(1-\theta_t)O(\norm{V})
\\
&+\left((1-\zeta_t)I-\frac{s_t}{\bar{\gamma}_t}D_{\vec{y}}\grad f(\vec{y})\right)\circ O(d(\vec{x},\vec{y}))+O(\norm{V})\circ O(d(\vec{x},\vec{y}))
\end{align}

\begin{align}
D_{\vec{v}}G&=D_{\vec{y}}G\circ D_{\vec{v}}\vec{y}+D_{\vec{v}}G\circ D_{\vec{v}}\vec{v}
\\
&=(1-\zeta_t)D_{\vec{v}}\vec{y}+\zeta_t I-\frac{s_t}{\bar{\gamma}_t}D_{\vec{y}}\grad f(\vec{y})D_{\vec{v}}\vec{y}+O(\norm{V})
\\
&=(1-\zeta_t)(\theta_t I+O(d(\vec{x},\vec{v})))+\zeta_t I-\frac{s_t}{\bar{\gamma}_t}D_{\vec{y}}\grad f(\vec{y})\circ (\theta_t I+O(d(\vec{x},\vec{v})))+O(\norm{V})
\\
&=(1-\zeta_t)\theta_t I+(1-\zeta_t)O(d(\vec{x},\vec{y}))+\zeta_t I-\frac{s_t}{\bar{\gamma}_t}(\theta_t D_{\vec{y}}\grad f(\vec{y})+D_{\vec{y}}\grad f(\vec{y})\circ O(d(\vec{x},\vec{v})))+O(\norm{V})
\\
&=((1-\zeta_t)\theta_t+\zeta_t)I-\frac{s_t\theta_t}{\bar{\gamma}_t}D_{\vec{y}}\grad f(\vec{y})+(1-\zeta_t)O(d(\vec{x},\vec{v}))-\frac{s_t}{\bar{\gamma}_t}D_{\vec{y}}\grad f(\vec{y})\circ O(d(\vec{x},\vec{v}))+O(\norm{V})
\end{align}

Suppose that $(\vec{x}^*,\vec{v}^*)=(\vec{x}^*,\vec{x}^*)$ is a critical point, then compute the Jacobian $D\psi(t,\vec{x},\vec{v})|_{(\vec{x}^*,\vec{v}^*)}$ at $(\vec{x}^*,\vec{v}^*)$ by letting $O(\cdot)=0$:

\begin{align}
D\psi(\vec{x}^*,\vec{v}^*)=
\left[
\begin{array}{cc}
(1-\theta_t)(I-\alpha_t D_{\vec{y}}\grad f(\vec{y}))&\theta_t(I-\alpha_t D_{\vec{y}}\grad f(\vec{y}))
\\
(1-\theta_t)((1-\zeta_t)I-\frac{s_t}{\bar{\gamma}_t}D_{\vec{y}}\grad f(\vec{y}))&((1-\zeta_t)\theta_t+\zeta_t)I-\frac{s_t\theta_t}{\bar{\gamma}_t}D_{\vec{y}}\grad f(\vec{y})
\end{array}
\right]_{\vec{y}=\vec{x}^*}
\end{align}

Next we compute the matrix $D_{\vec{y}}\grad f(\vec{y})$ and evaluate it at saddle point. We firstly recall that the Riemannian gradient $\grad f(\vec{y})$ in local coordinate system is
\[
\grad f(\vec{y})=\left(g^{1j}\frac{\partial f}{\partial y_j},...,g^{dj}\frac{\partial f}{\partial y_j}\right)=(g^{ij})\cdot\nabla f(\vec{y})
\]
where 
\[
g^{ij}\frac{\partial f}{\partial y_j}=\sum_{j=1}^dg^{ij}\frac{\partial f}{\partial y_j}
\]
according to Einstein's convention and $(g^{ij})$ denote the inverse of the metric matrix $(g_{ij})$. The differential of the Riemannian gradient in local coordinates is computed as follows,
\begin{align}
D_{\vec{y}}\grad f(\vec{y})&=D_{\vec{y}}((g^{ij})\cdot\nabla f(\vec{y}))
\\
&=\left[
\begin{array}{ccc}
\frac{\partial g^{1j}}{\partial y_1}\frac{\partial f}{\partial y_j}&\dotsm&\frac{\partial g^{1j}}{\partial y_d}\frac{\partial f}{\partial y_j}
\\
\vdots&&\vdots
\\
\frac{\partial g^{dj}}{\partial y_1}\frac{\partial f}{\partial y_j}&\dotsm&\frac{\partial g^{dj}}{\partial y_d}\frac{\partial f}{\partial y_j}
\end{array}
\right]
+
\left[
\begin{array}{ccc}
g^{1j}\frac{\partial^2f}{\partial y_1\partial y_j}&\dotsm&g^{1j}\frac{\partial^2f}{\partial y_d\partial y_j}
\\
\vdots&&\vdots
\\
g^{dj}\frac{\partial^2f}{\partial y_1\partial y_j}&\dotsm&g^{dj}\frac{\partial^2f}{\partial y_d\partial y_j}
\end{array}
\right].
\end{align}
Since at a saddle point, or more generally, a critical point, the partial derivatives $\frac{\partial f}{\partial y_j}$ are all equal to zero, under any local coordinate systems, then we have
\[
D_{\vec{y}}\grad f(\vec{y})=\left[
\begin{array}{ccc}
g^{1j}\frac{\partial^2f}{\partial y_1\partial y_j}&\dotsm&g^{1j}\frac{\partial^2f}{\partial y_d\partial y_j}
\\
\vdots&&\vdots
\\
g^{dj}\frac{\partial^2f}{\partial y_1\partial y_j}&\dotsm&g^{dj}\frac{\partial^2f}{\partial y_d\partial y_j}
\end{array}
\right]
\]
On the other hand, by choosing local coordinate systems properly, e.g. in a normal neighborhood at the saddle point $\vec{x}^*$, we can make $g_{ij}=g^{ij}=\delta_{ij}$. So under such special coordinate system, the differential of gradient is identical to the Euclidean Hessian at the critical point, i.e.
\[
D_{\vec{y}}\grad f(\vec{y})=\nabla^2f(\vec{y})
\]
and then
\[
D\psi(t,\vec{x}^*,\vec{v}^*)=\left[
\begin{array}{cc}
(1-\theta_t)(I-\alpha_t\nabla^2 f(\vec{y}))&\theta_t(I-\alpha_t \nabla^2f(\vec{y}))
\\
(1-\theta_t)((1-\zeta_t)I-\frac{s_t}{\bar{\gamma}_t}\nabla^2f(\vec{y}))&((1-\zeta_t)\theta_t+\zeta_t)I-\frac{s_t\theta_t}{\bar{\gamma}_t}\nabla^2f(\vec{y})
\end{array}
\right]_{\vec{y}=\vec{x}^*}
\]
in the normal neighborhood of $\vec{x}^*$. Since the Hessian is diagonalizable, the block matrices of $D\psi(t,\vec{x}^*,\vec{v}^*)$ are simultaneously diagnolizable, the rest proof is a consequence of following two lemmas.
\end{proof}

\begin{lemma}\label{Jacobian:diag}%%%%%%%% diagonalizable
$D\psi(t,\vec{x}^*,\vec{v}^*)$ is diagonalizable.
\end{lemma}

\begin{proof}
Under similar transformation 

\[
\tilde{D}:=
\left[
\begin{array}{cc}
C^{-1}&0
\\
0&C^{-1}
\end{array}
\right]D\psi(t,\vec{x}^*,\vec{v}^*)
\left[
\begin{array}{cc}
C&0
\\
0&C
\end{array}
\right]
\]

 $D\psi(t,\vec{x}^*,\vec{v}^*)$ is similar to 
\[
\tilde{D}=
\left[
\begin{array}{cc}
(1-\theta_t)(I-\alpha_t H)&\theta_t(I-\alpha_t H)
\\
(1-\theta_t)\left((1-\zeta_t)I-\frac{s_t}{\bar{\gamma}_t}H\right)&\left((1-\zeta_t)\theta_t+\zeta_t\right)I-\frac{s_t\theta_t}{\bar{\gamma}_t}H
\end{array}
\right]
\]
where $H=\diag\{\lambda_i\}$. The determinant of $\tilde{D}-xI_{2d}$ is
\begin{align}
&\det(\tilde{D}-xI)
\\
&=\det\left(
\begin{array}{cc}
(1-\theta_t)(I-\alpha_t H)-xI&\theta_t(I-\alpha_t H)
\\
(1-\theta_t)\left((1-\zeta_t)I-\frac{s_t}{\bar{\gamma}_t}H\right)&\left((1-\zeta_t)\theta_t+\zeta_t\right)I-\frac{s_t\theta_t}{\bar{\gamma}_t}H-xI
\end{array}
\right)
\\
&=\det\left(\left((1-\theta_t)(I-\alpha_t H)-xI\right)\left(((1-\zeta_t)\theta_t+\zeta_t)I-\frac{s_t\theta_t}{\bar{\gamma}_t}H-xI\right)\right.
\\
&\left.-\theta_t(1-\theta_t)(I-\alpha_t H)\left((1-\zeta_t)I-\frac{s_t}{\bar{\gamma}_t}H\right)\right)
\\
&=\det\left(\diag\{\left[(1-\theta_t)(1-\alpha_t\lambda_i)-x\right]\left[\left((1-\zeta_t)\theta_t+\zeta_t-\frac{s_t\theta_t}{\bar{\gamma}_t}\lambda_i\right)-x\right]\right.
\\
&\left.-\theta_t(1-\theta_t)(1-\alpha_t\lambda_i)\left((1-\zeta_t)-\frac{s_t}{\bar{\gamma}_t}\lambda_i\right)\}\right)
\end{align}
The entry of diagonal matrix equals
\begin{align}
&x^2-\left[(1-\theta_t)(1-\alpha_t\lambda_i)+\left((1-\zeta_t)\theta_t+\zeta_t-\frac{s_t\theta_t}{\bar{\gamma}_t}\lambda_i\right)\right]x
\\
&+(1-\theta_t)(1-\alpha\lambda_i)\left((1-\zeta_t)\theta_t+\zeta_t-\frac{s_t\theta_t}{\bar{\gamma}_t}\lambda_i\right)
\\
&-\theta_t(1-\theta_t)(1-\alpha_t\lambda_i)\left((1-\zeta_t)-\frac{s_t}{\bar{\gamma}_t}\lambda_i\right)
\\
&=x^2-\left[(1-\theta_t)(1-\alpha_t\lambda_i)+\left((1-\zeta_t)\theta_t+\zeta_t-\frac{s_t\theta_t}{\bar{\gamma}_t}\lambda_i\right)\right]x
\\
&+(1-\theta_t)(1-\alpha_t\lambda_i)\left[(1-\zeta_t)\theta_t+\zeta_t-\frac{s_t\theta_t}{\bar{\gamma}_t}\lambda_i-\theta_t(1-\zeta_t)+\frac{s_t\theta_t}{\bar{\gamma}_t}\lambda_i\right]
\\
&=x^2-\left[(1-\theta_t)(1-\alpha_t\lambda_i)+\left((1-\zeta_t)\theta_t+\zeta_t-\frac{s_t\theta_t}{\bar{\gamma}_t}\lambda_i\right)\right]x+(1-\theta_t)(1-\alpha_t\lambda_i)\zeta_t
\end{align}
The discriminant of above quadratic function is
\[
\Delta=\left[(1-\theta_t)(1-\alpha_t\lambda_i)+\left((1-\zeta_t)\theta_t+\zeta_t-\frac{s_t\theta_t}{\bar{\gamma}_t}\lambda_i\right)\right]^2-4(1-\theta_t)(1-\alpha_t\lambda_i)\zeta_t.
\]
Notice that when $\alpha_t=0$, we have $\Delta=(1+\zeta_t)^2-4\zeta_t=(\zeta_t-1)^2$, since $s_t,\theta_t\rightarrow 0$ as $\alpha_t\rightarrow 0$. On the other hand, since $\zeta_t=\frac{(1-s_t)\gamma_t}{\bar{\gamma}_t}=\frac{1-s_t}{1+\beta_t}\ne 1$, then $(\zeta_t-1)^2>1$ always holds. Therefore, there exists small interval $I$ to which $\alpha_t$ belongs such that the quadratic polynomial has two real roots, and this implies that the polynomial
\[
\det\left(\tilde{D}-xI_{2d}\right)
\]
has $2d$ real roots. Thus the matrix $\tilde{D}$ is diagonalizable.
\end{proof}

\begin{lemma}\label{DD}
 $D\psi(t,\vec{x}^*,\vec{v}^*)$ has an eigenvalue greater than 1 when the step size $\alpha_t$ is taken to be small enough..
 \end{lemma}
 
 \begin{proof}
Since $\nabla^2 f(\vec{x}^*)$ is a real symmetric matrix, it can be diagonalized by a similar transformation. Assume we have 
$$C \cdot \nabla^2 f(\vec{x}^*) \cdot C^{-1} = \lambda_f $$
where $\lambda_f$ is 
\[\left(\begin{array}{cccc}
   \lambda_1 &   &   &   \\
   &\lambda_2 &  &  \\
   &  & \ddots &  \\
  & & &\lambda_d
\end{array}\right)\]
and $\{\lambda_i\}_{i \in [n]}$ are eigenvalues of $\nabla^{2} f(\vec{x}^{\ast})$, and we also assume $\lambda_d < 0$. 

Then do a similar transformation to $D\psi(t,\vec{x}^*,\vec{v}^*)$ by matrix 
$$\left[
\begin{array}{cc}
C&0
\\
0&C
\end{array}
\right]
$$
we get
$$ 
\left[
\begin{array}{cc}
C&0
\\
0&C
\end{array}
\right]
\cdot
D\psi(t,\vec{x}^*,\vec{v}^*)
\cdot
\left[
\begin{array}{cc}
C^{-1}&0
\\
0&C^{-1}
\end{array}
\right]
= 
\left[
\begin{array}{cc}
(1-\theta_t) (I_{d \times d} - \alpha_t \lambda_f) & \theta_t(I_{d \times d}-\alpha_t \lambda_f) \\ (1-\theta_t)((1-\zeta_t)I_{d \times d} - \frac{s}{\bar{\gamma}_t} \lambda_f) & ((1-\zeta_t)\theta_t + \zeta_t)I_{d \times d} - \frac{s \theta_t}{\bar{\gamma}_t} \lambda_f 
\end{array}
\right]
$$
We denote the above matrix by $\widetilde{D\psi(t,\vec{x}^*,\vec{v}^*)}$.

Then we calculate the determinate of $$xI_{2d \times 2d} - \widetilde{D\psi(t,\vec{x}^*,\vec{v}^*)} $$
 which equals to $\prod^d_{i = 1} g_i(x)$, and 
$$g_i(x) = x^2 + b_i x + c_i$$
where $b_i = (\alpha_t(1-\theta_t) + \frac{s_t \theta_t}{\bar{\gamma}_t}) \lambda_i -\zeta_t(1-\theta_t) -1$ and $c_i = \zeta_t(1-\theta_t) (1-\alpha_t \lambda_i)$.
Recall that we assume that $\lambda_d<0$, we focus on the quadratic polynomial $g_d(x)$ and analyze its roots. For 
\[
g_d(x)=x^2+b_dx+c_d
\]
where 
\[
b_d=(\alpha_t(1-\theta_t)+\frac{s_t\theta_t}{\bar{\gamma}_t})\lambda_d-\zeta_t(1-\theta_t)-1
\]
and 
\[
c_d=\zeta_t(1-\theta_t)(1-\alpha_t\lambda_d).
\]
%the axis of symmetry of the parabola $g_d(x)$ is 
%\[
%\frac{x_1+x_2}{2}=-\frac{b_d}{2}=-\frac{1}{2}\left((\alpha(1-\theta)+\frac{s\theta}{\bar{\gamma}})\lambda_d-\zeta(1-\theta)-1\right),
%\]
%when $\alpha\rightarrow 0$, we have that $-\frac{b_d}{2}\rightarrow 1$. 

By the proof of Lemma \ref{Jacobian:diag}, we know that when the step size is small enough, all $g_i(x)$ have two real roots. Thus the larger root of $g_d(x)$, i.e.
\[
\frac{-b_d+\sqrt{b_d^2-4c_d}}{2}
\]
is greater than 1 if the following holds:
\[
-b_d(c_d+1)+(c_d+1)^2>-b_d^2c_d.
\]
We note that it trivially holds since $b_d\rightarrow -2$ and $c_d\rightarrow 1$ as $\alpha_t\rightarrow 0$, since $\alpha_t(1-\theta_t)+\frac{s\theta_t}{\bar{\gamma}_t}\rightarrow 0$, $\zeta_t\rightarrow 1$ and $\theta_t\rightarrow 0$. So the proof completes.
%Recall that we have assumed that $\lambda_d < 0$, and now we claim that $g_d(x)$ has a root greater than 1. That is because under the condition that $\alpha - \frac{s}{\bar{\gamma}} > 0$, it's easily to check that the discriminate of the quadratic polynomial $g_d(x)$ is greater then $0$, and by Vieta's theorem, the summation of two roots of $g_d(x)$ equals to $-b_d$, which is always greater than 1.   
\end{proof}

So the eigenvalue of magnitude greater than 1 in $I-\alpha_t\nabla^2 f(\vec{x}^*)$ is also an eigenvalue of $D\psi(t,\vec{x}^*,\vec{v}^*)$, which is a consquence of the assumption that $\vec{x}^*$ is saddle point. The proof of Lemma \ref{lemma:eigenvalue} completes.

\twocolumn

%\bibliography{ijcai22}

\end{document}

%% file: intronew.tex
\section{Introduction}

In this paper we consider non-convex optimization problem with constraint that is a product of simplices, i.e.,
\begin{align}
\min_{\vec{x}\in \Delta_1\times...\times \Delta_n} f(\vec{x})
\end{align}
where $f:\Delta_1\times...\times \Delta_n\rightarrow\mathbb{R}$ is a sufficiently smooth function and
\[\Delta_i=\left\{(x_{i1},...,x_{id}):\sum_{s=1}^dx_{is}=1,x_{is}\ge 0\right\},\]
%where $\mathcal{M}=\Delta_1\times...\times \Delta_n$ and $f:\mathcal{M}\rightarrow\mathbb{R}$ is a sufficiently smooth function. Especially, we focus on the optimization with Riemannian Accelerated Gradient Descent and its saddle avoidance behavior, with applications in Multiplicative Weights Update (MWU) \cite{AHK12}, which is commonly used in various fields including game theory, optimization, machine learning and multi-agent systems.
Especially, we are interested in the Multiplicative Weights Update (MWU) algorithm ~\cite{AHK12}, which is commonly used in various fields including game theory, optimization, machine learning and multi-agent systems. We propose an Accelerated Multiplicative Weights Update (A-MWU), and study its saddle point avoidance behavior from a general perspective of Riemannian Accelerated Gradient Descent \cite{ZS18}.

%There has been a growing interest in optimization on Riemannian manifolds due to the generic existence of non-Euclidean structure in Machine Learning problems. For example, matrix manifolds play a fundamental role in computer vision and pattern recognition. Recently, Geometric Deep Learning ~\cite{LeCun17} as the intersection of differential geometry and deep learning has led to new challenges of learning embeddings for non-Euclidean data, modeling hierarchical structure, and optimization on Riemannian manifolds \cite{NickelKiela17,NickelKiela18}.

 Escaping saddle points in non-convex optimization has been studied extensively by the Machine Learning community \cite{GHJY2015,Jin17,JNJ18,LSJR16,LP19,panageas2016gradient,CB19,SFF19,PPW19a,SLQ19}. The first-order optimization algorithms can be studied from a dynamical system perspective and results in\cite{LSJR16,LP19} guarantee that the algorithms asymptotically avoid saddle points in probability 1 with random initialization. The main technique of using Center-Stable Manifold theorem was extended to heavy-ball algorithm on Euclidean space by \cite{SLQ19}. 
 
 Despite the aforementioned progresses, the saddle avoidance of accelerated algorithms in non-Euclidean setting is less studied, especially for the case where the step-size is varying with time. The only result on variable step-size is given by \cite{PPW19a}, but the Riemannian Accelerated Gradient Gescent is omitted. It is known that MWU on the simplex is a special case of mirror descent with entropy regularizer, and the mirror descent algorithm has many applications in optimization \cite{DGSX2011,JN2011}.
 Using this connection, an accelerated version of mirror descent is provided by \cite{KBB}, in which the authors discretize a system of ODEs of continuous-time mirror descent to get acceleration. However, \cite{KBB} focuses on convex optimization only and leaves the saddle avoidance in non-convex optimization.%However, the experiments in \cite{KBB} indicate that the algorithm may oscillate.  
\  \cite{ICMLPPW19} proved that MWU almost always converges to second-order stationary points with constant step-size, but the result for accelerated MWU with variable step-size is missing in the literatures.

%\onecolumn
\begin{table*}
\centering
\begin{tabular}{ |p{3.5cm}||p{2cm}|p{3cm}|p{3cm}|p{3cm}|  }
 \hline
 \multicolumn{5}{|c|}{Table 1: Comparison to related results} \\
 \hline
 & Acceleration &Saddle Avoidance& Variable Step-size& Types of Manifolds\\
 \hline
Accelerated MD \cite{KBB}   & \cmark    &\xmark& \cmark  & Convex set\\
 Heavy-ball  \cite{SLQ19}&  \cmark   & \cmark  &\xmark& Euclidean space\\
MWU  \cite{ICMLPPW19} &\xmark & \cmark&  \xmark& Simplices\\
 RAGD  \cite{ZS18}   &\cmark & \xmark&  \cmark& General\\
 AMWU (this work)&   \cmark  & \cmark&\cmark& General, Simplices\\
 \hline
\end{tabular}
\end{table*}
%\twocolumn
Motivated by the question that if there exists an accelerated version of MWU which also provably avoids saddle point, we investigate the Riemannian Accelerated Gradient Descent (RAGD) proposed in \cite{ZS18}. The main results and contributions of this paper are the following:
\begin{itemize}
\item We propose an Accelerated Multiplicative Weights Update which is derived and simplified from RAGD of \cite{ZS18};
\item We prove that the RAGD of \cite{ZS18} avoids saddle points and moreover, this provides the first saddle avoidance result of Accelerated Multiplicative Weights Update.
\end{itemize} %and show that their algorithm provably avoids saddle points of locally geodesically convex functions on Riemannian manifold. Furthermore, we propose an accelerated version of Multiplicative Weights Update based on the Riemannian geometry of the simplex, and apply the saddle avoidance results to Accelerated MWU. 
%Our contributions compared to the most relevant results in literature are illustrated in Table 1.

\begin{comment}
\paragraph{More related works.}We mention the other line of work by \cite{GHJY2015,Jin17}, the authors show a non-asymptotic guarentee for the \emph{perturbed gradient descent} escape from saddle points efficiently. \cite{ge2015escaping} present the first polynomial guarantees for a perturbed version of gradient descent. Based on this result,  \cite{Jin17} proved that perturbed gradient descent can find a $\epsilon$-second-order stationary point in $\tilde{O}(1/\epsilon^2)$, their proof relies on a new characterization of the geometry around saddle points. \cite{JNJ18} proved perturbed accelerated  gradient descent can escape saddle points faster than gradient descent, and find $\epsilon$-second-order stationary point in $\tilde{O}(1/\epsilon^\frac{7}{4})$. Recently their result has been generalized to first-order method on Riemannian manifold by \cite{Boumal2020}. The result of \cite{Jin17} on perturbed gradient descent is generalized to Riemannian manifold independently by \cite{CB19} and \cite{SFF19}. The stochastic accelerated mirror descent is given by Xu et al. in \cite{XWG18}, based on a continuous-time SDE argument.
\end{comment}

%% file: mx.tex
\newpage
\section{Preliminaries}
\paragraph{Riemannian metric and geodesic.} A $d$-dimension Riemannian manifold $(M,g)$ is real, smooth $d$-dimension manifold $M$ equipped with a Riemannian metric $g$. For each $\vec{x}\in M$, let $T_{\vec{x}}M$ denote the tangent space at $\vec{x}$. The metric $g$ induces a inner product $\langle\cdot,\cdot\rangle_{\vec{x}}:T_{\vec{x}}M\times T_{\vec{x}}M\rightarrow\mathbb{R}$. We call a curve $\gamma(t):[0,1]\rightarrow M$ a geodesic if it satisfies 
\begin{itemize}
\item The curve $\gamma(t)$ is parametrized with constant speed, i.e. $\norm{\frac{d}{dt}\gamma(t)}_{\gamma(t)}$ is constant for $t\in[0,1]$.
\item The curve is locally length minimized between $\gamma(0)$ and $\gamma(1)$.
\end{itemize}

\paragraph{Exponential and logarithmic map.}The exponential map $\Exp_{\vec{x}}(\vec{v})$ maps $\vec{v}\in T_{\vec{x}}M$ to $\vec{y}\in M$ such that there exists a geodesic $\gamma$ with $\gamma(0)=\vec{x}$, $\gamma(1)=\vec{y}$ and $\gamma'(0)=\vec{v}$. For $\vec{x}\in M$, let $\Log_{\vec{x}}$ denote the logarithmic map at $\vec{x}$,
\[
\Log_{\vec{x}}(\vec{y})=\argmin_{\vec{u}\in T_{\vec{x}}M} \ \ \text{subject to}\ \ \Exp_{\vec{x}}(\vec{u})=\vec{y},
\]
with domain such that this is uniquely defined.

\begin{remark}
The existence and uniqueness of geodesic is guaranteed by the Fundametal Theorem of ODE, one cannot expect an explicit expression of the geodesic, usually we need to apply a numerical scheme to approximate the solution of the geodesic equation. However, there are cases for which the geodesic equation can be expressed with a closed form, we give the following examples.
\end{remark}

\begin{comment}
\begin{example}
Let $S^n=\{\vec{x}\in\mathbb{R}^{n+1}:\norm{\vec{x}}=1\}$ be the unit sphere embedded in Euchlidean space. The exponential map $\Exp_{\vec{x}}(\vec{v}):T_{\vec{x}}M\rightarrow M$ is given by 
\[
\Exp_{\vec{x}}(\vec{v})=\cos(\norm{\vec{v}})\vec{x}+\sin(\vec{\norm{\vec{v}}})\frac{\vec{v}}{\norm{\vec{v}}}
\]
where $\norm{\vec{v}}=\sqrt{\langle\vec{v},\vec{v}\rangle}$ is computed from the Euclidean inner product. 
\end{example}
\end{comment}

\paragraph{Geodesically convex.}A set $U\subset M$ is geodesically convex if for any $\vec{x},\vec{y}\in U$, there is a geodesic $\gamma$ with $\gamma(0)=\vec{x}$, $\gamma(1)=\vec{y}$ and $\gamma(t)\in U$ for $t\in[0,1]$.

Let $U$ be a geodesically convex subset of $M$. A function $f: M\rightarrow\mathbb{R}$ is called geodesically convex on $U$ if for any $\vec{x},\vec{y}\in M$ and any geodesic $\gamma$ such that $\gamma(0)=\vec{x}$, $\gamma(1)=\vec{y}$ and $\gamma(t)\in U$ for all $t\in[0,1]$, it holds that
\[
f(\gamma(t))\le (1-t)f(\vec{x})+tf(\vec{y}).
\]
A function $f: M\rightarrow\mathbb{R}$ is called geodesically $\mu$-convex on $U$ if for any $\vec{x},\vec{y}\in U$ and gradient $\grad f(\vec{x})$ at $\vec{x}$, it holds that
\begin{equation}\label{geodconv}
f(\vec{y})\ge f(\vec{x})+\langle\grad f(\vec{x}),\Log_{\vec{x}}(\vec{y})\rangle+\frac{\mu}{2}\norm{\Log_{\vec{x}}(\vec{y})}^2,
\end{equation}
where we assume the logarithmic map is well defined in $U$.

\paragraph{Riemannian Gradient and Hessian.}
 For differentiable function $f:M\rightarrow\BR$, $\grad f(\vec{x})\in T_{\vec{x}}M$ denotes the Riemannian gradient of $f$ that satisfies $\frac{d}{dt}f(\gamma(t))=\langle\gamma'(t),\grad f(\vec{x})\rangle$ for any differentiable curve $\gamma(t)$ passing through $\vec{x}$. The local coordinate expression of gradient is useful in our analysis.
\begin{equation}\label{grad}
\grad f(\vec{x})=\left(\sum_jg^{1j}(\vec{x})\frac{\partial f}{\partial x_j},...,\sum_jg^{dj}(\vec{x})\frac{\partial f}{\partial x_j}\right)
\end{equation}
where $g^{ij}(\vec{x})$ is the $ij$-th entry of the inverse of the metric matrix $\{g_{ij}(\vec{x})\}$ at each point. 

The Hessian of $f$ is the covariant derivative of the gradient vector field: $\Hess f(\vec{x})[\vec{u}]=\nabla_{\vec{u}}\grad f(\vec{x})$ for any vector field $\vec{u}$ on $M$.

\paragraph{Strict saddle point.} A strict saddle point $\vec{x}^*$ of $f:M\rightarrow\mathbb{R}$ satisfies
\[
\norm{\grad f(\vec{x}^*)}=0 \ \ \text{and} \ \ \lambda_{\min}(\Hess f(\vec{x}^*))<0.
\]

\paragraph{Retraction.}
A retraction on a manifold $M$ is a smooth mapping $\Retr$ from the tangent bundle $TM$ to $M$ satisfying properties 1 and 2 below: Let $\Retr_{\vec{x}}:T_{\vec{x}}M\rightarrow M$ denote the restriction of $Retr$ to $T_{\vec{x}}M$.
\begin{enumerate}
\item $\Retr_{\vec{x}}(0)=\vec{x}$, where $0$ is the zero vector in $T_{\vec{x}}M$.
\item The differential of $\Retr_{\vec{x}}$ at $0$ is the identity map.
\end{enumerate}
Then the Riemannian gradient descent with stepsize $\alpha $ is given as
\begin{equation}\label{GD:Riemannian}
\vec{x}_{t+1}=\Retr_{\vec{x}_t}(-\alpha\grad f(\vec{x}_t)).
\end{equation}

\paragraph{Multiplicative Weights Update.} The paper focuses on Accelerated Multiplicative Weights Update, for completeness, we recall the linear variant of MWU. Suppose that $\vec{x}_i=(x_{i1},...,x_{id_i})$ is in the $i$-th component of $\Delta_1\times...\times\Delta_n$. Assume that $\vec{x}(t)$ is the $t$-th iterate of MWU, the algorithm is written as follows:
\begin{equation}\label{MWUclassic}
x_{ij}(t+1)=x_{ij}(t)\frac{1-\alpha\frac{\partial f}{\partial x_{ij}}}{1-\alpha\sum_{s}x_{is}(t)\frac{\partial f}{\partial x_{is}}},
\end{equation}
where $j\in\{1,...,d_i\}$.

\section{Algorithm}

We give the Single-Agent version of the Accelerated Multiplicative Weights Update (A-MWU) in Algorithm \ref{alg:C}, which is derived from the Riemannian Accelerated Gradient Descent based on the geometry of the positive orthant $\mathbb{R}^d_+$. %The multi-agent A-MWU is implemented as follows: at the $t$'th step, algorithm \ref{alg:C} is executed for each $i\in[n]$. 
We leave the full multi-agent A-MWU algorithm to Appendix.

 \begin{algorithm}
\caption{Single-Agent A-MWU }
\label{alg:C}
\begin{algorithmic}
\STATE {input : $\vec{x}_0, \vec{v}_0, 0<c \le\alpha_t< \frac{1}{L}, \beta_t >0$,$\delta >0$}, 
\\
\REPEAT
\STATE Compute $s_t\in(0,1)$ from the equation $s_t^2=\alpha_t((1-s_t)\gamma_t+s_t\mu)$,
\\
Set $\bar{\gamma}_{t+1}=(1-s_t)\gamma_t+s_t\mu$, $\gamma_{t+1}=\frac{1}{1+\beta_t}\bar{\gamma}_{t+1}$,  
\\$S=\left(\prod_{i=1}^d\frac{x_i(t)}{y_i(t)}\right)^{1/d}$
\\
Set $y_i(t+1)=\frac{x_i(t)\exp\left(\frac{s_t\gamma_t}{\gamma_t+s_t\mu}\ln\left(S\frac{v_i(t)}{x_i(t)}\right)\right)}{\sum_jx_j(t)\exp\left(\frac{s_t\gamma_t}{\gamma_t+s_t\mu}\ln\left(S\frac{v_j(t)}{x_j(t)}\right)\right)}$
%Set $\vec{y}_{t}=\Exp_{\vec{x}_t}\left(\frac{s_t\gamma_t}{\gamma_t+s_t\mu}\Log_{\vec{x}_t}(\vec{v}_t)\right)$
\\
Set $x_i(t+1)=y_i(t)\frac{1-\alpha_t\frac{\partial f}{\partial x_i}(\vec{y}_t)}{1-\alpha_t\sum_j\frac{\partial f}{\partial x_j}(\vec{y}_t)}$
%Set $\vec{x}_{t+1}=\Exp_{\vec{y}_t}(-\alpha_t\grad f(\vec{y}_t))$
\\
Compute $S'=\left(\prod_{i=1}^d\frac{y_i(t)}{v_i(t)}\right)^{1/d}$, 
 \[
 \begin{split}u_i&=\frac{(1-s_t)\gamma_t}{\bar{\gamma}_t}\ln\left(S'\frac{v_i(t)}{y_i(t)}\right)\\
 &\quad+y_i(t)\frac{1-\alpha_t\frac{\partial f}{\partial x_i}(\vec{y}_t)}{1-\alpha_t\sum_jy_j(t)\frac{\partial f}{\partial x_j}(\vec{y}_t)}-y_i(t)
 \end{split}
 \]
Set $v_i(t+1)=\frac{y_i(t)\exp(u_i)}{\sum_jy_j(t)\exp(u_j)}$
%Set $\vec{v}_{t+1}=\Exp_{\vec{y}_t}\left(\frac{(1-s_t)\gamma_t}{\bar{\gamma}_t}\Log_{\vec{y}_t}(\vec{v}_t)-\frac{s_t}{\bar{\gamma}_t}\grad f(\vec{y}_t)\right)$ 

\UNTIL{$\norm{\grad f(\vec{y}_t)} \le \delta$}
\end{algorithmic}
\end{algorithm}
%\end{comment}
Even the derivation of Algorithm \ref{alg:C} is based on the Riemannian Accelerated Gradient Descent, the above implementation is simplified compared to the original algorithm in \cite{ZS18}. Since it is impossible to have an explicit form of exponential and logarithmic map in general, all the updates of the algorithm are computed based on the property of Riemannian geometry of the positive orthant $\mathbb{R}^d_+=\{\vec{x}:x_i>0 \ \ \text{for all}\ \ i\in[d]\}$, \cite{Shahshahani,HofSig}. We recall the background that is necessary for the interpretation that MWU is indeed a manifold gradient descent. Formally, the positive orthant $\mathbb{R}^d_+$ is endowed with a Riemannian metric (called the \textbf{Shahshahani metric}) whose metric matrix $\{g_{ij}(\vec{x})\}$ is diagonal with $g_{ii}(\vec{x})=\frac{\abs{\vec{x}}}{x_i}$ where $\abs{\vec{x}}=\sum_jx_j$, %or equivalently,
%\[
%g(\vec{x})=\left(
%\begin{array}{ccc}
%\frac{\abs{\vec{x}}}{x_1}&&\text{\Large 0}
%\\
%&\ddots&
%\\
%\text{\Large 0}&&\frac{\abs{\vec{x}}}{x_d}
%\end{array}\right)
%\]
 %Usually a manifold is defined by a collection of coordinate charts, here $\mathbb{R}^d_+$ is considered a single chart manifold with a non-Euclidean structure, and 
 The positive orthant $\mathbb{R}^d_+$ with the Shahshahani metric is call a Shahshahani manifold. The tangent spaces $T_{\vec{x}}\mathbb{R}^d_+$ for all $\vec{x}\in \mathbb{R}^d_+$ are all identified with $\mathbb{R}^d$. Now consider a differentiable function $f:\mathbb{R}^d_+\rightarrow\mathbb{R}$, one can define the gradient of $f$ at each point with respect to the Shahshahani metric. The following local-coordinate expression of the Shahshahani gradient for each $\vec{x}\in\mathbb{R}^d_+$ is straightforward from (\ref{grad}):
\[
\grad f(\vec{x})=g^{-1}\cdot\nabla f(\vec{x})=\left(\frac{x_1}{\abs{\vec{x}}}\frac{\partial f}{\partial x_1},...,\frac{x_d}{\abs{\vec{x}}}\frac{\partial f}{\partial x_d}\right).
\]
 We know from \cite{HofSig} that the exponential map on the positive simplex $\Delta_+\subset\mathbb{R}^n_+$ is
\[
\Exp_{\vec{x}}(\vec{v})=\left(\frac{x_1e^{v_1}}{\sum x_ie^{v_i}},...,\frac{x_ne^{v_n}}{\sum x_ie^{v_i}}\right),
\]
where $\vec{v}=(v_1,...,v_n)$ is a tangent vector of $\Delta_+$ at point $\vec{x}$.

Moreover, the Multiplicative Weights Update given by (\ref{MWUclassic}) can be understood as a retraction on $\Delta_+$ and the logarithmic map can also be computed explicitly. These makes the formulation of Algorithm \ref{alg:C} possible, and we leave the details in Appendix.

%Note that although many Riemann geometry based optimization algorithms contain calculations of  $\Exp$ and $\Log$ maps, which are difficult in practice. But here, due to the special geometry structure of  $\mathbb{R}^d_+$, the $\Exp$ and $\Log$ map can be written down explicitly, thus in  Algorithm \ref{alg:C} every step can be executed directly. The algorithm is based on the Accelerated Riemannian Gradient Descent in \cite{ZS18},  we leave the description of geometry on $\mathbb{R}^d_+$ and the general Accelerated Riemannian Gradient Descent  algorithm to the next section.

\begin{comment}
The theoretical guarantee in saddle avoidance of Riemannian Accelerated Gradient Descent implies that the Accelerated Multiplicative Weights Update avoids saddle points almost always, i.e., the following corollary is immediate from Theorem \ref{SMT:RAGD2} and Corollary \ref{distributed}.

\begin{corollary}\label{AMWU}
Let $\mathcal{M}=\Delta_1\times...\times\Delta_n$, where 
\[
\Delta_i=\left\{\vec{x}_i\in\mathbb{R}^{d_i}:\sum_{s=1}^{d_i}x_{is}=1,x_{is}\ge0\right\}.
\]
Suppose $f:\mathcal{M}\rightarrow\mathbb{R}$ is $C^2$ and geodesically convex w.r.t. the product Shahshahani metric. Then the Accelerated Multiplicative Weights Update algorithm avoids interior saddle points almost always, i.e., randomly choose an initial point, the algorithm will avoid interior saddle points with probability one. 
\end{corollary}
\end{comment}

\section{Saddle Avoidance Analysis}
%\subsection{Accelerated Gradient Descent on Product Manifold}
We show that the Accelerated Gradient Descent on the product manifold $M_1\times...\times M_n$ actually avoids saddle points almost always. The algorithm is extended from the Accelerated Riemannian Gradient Descent of \cite{ZS18} with the structure of product manifold. %Before going through the detailed analysis, we recall basic facts on the Riemannian geometry of the positive orthant.
%\paragraph{Riemannian geometry of the positive orthant.}

Suppose that each manifold component $M_i$ is equipped with a Riemannian metric $g_i$, then the product manifold $M_1\times...\times M_n$ has the product metric $g=g_1\otimes...\otimes g_n$ whose metric matrix is blocked matrix with $g_i$'s the non-trivial blocks. Then the gradient is also the Cartesian product of the gradients of each component manifold $M_i$, i.e. 
\[
\grad f(\vec{x})=(\grad_{\vec{x}_1} f(\vec{x}),...,\grad_{\vec{x}_n} f(\vec{x})).
\]

Recall the Accelerated Riemannian Gradient Descent in \cite{ZS18} as follows, 
\begin{align}\label{RAGDsra}
\vec{y}_{t}&=\Exp_{\vec{x}_t}\left(\frac{s_t\gamma_t}{\gamma_t+s_t\mu}\Log_{\vec{x}_t}(\vec{v}_t)\right)
\\
\vec{x}_{t+1}&=\Exp_{\vec{y}_t}(-\alpha_t\grad f(\vec{y}_t))
\\
\vec{v}_{t+1}&=\Exp_{\vec{y}_t}\left(\frac{(1-s_t)\gamma_t}{\bar{\gamma}_t}\Log_{\vec{y}_t}(\vec{v}_t)-\frac{s_t}{\bar{\gamma}_t}\grad f(\vec{y}_t)\right)
\end{align}
where the $0<c \le\alpha_t< \frac{1}{L}$, $\beta_t >0$, and $s_t$, $\gamma_t$, $\bar{\gamma}_t$ are computed according to $s_t\in(0,1)$, such that $s_t^2=\alpha_t((1-s_t)\gamma_t+s_t\mu)$, and 
\[
\bar{\gamma}_{t+1}=(1-s_t)\gamma_t+s_t\mu, \ \ \ \gamma_{t+1}=\frac{1}{1+\beta_t}\bar{\gamma}_{t+1}.
\]

In order to write the distributed Accelerated Riemannian Gradient Descent, we simply compute the exponential and logarithmic map component-wise. We denote 
\[
\overrightarrow{\vec{x}}=(\vec{x}_1,...,\vec{x}_n)\in M_1\times,...,\times M_n
\]
and the $t$'th update
\[
\overrightarrow{\vec{x}}_t=(\vec{x}_{1,t},...,\vec{x}_{n,t}), \ \ \overrightarrow{\vec{v}}_t=(\vec{v}_{1,t},...,\vec{v}_{n,t})
\]
are computed according to (\ref{RAGDsra}) in a distributed manner.

\begin{comment}
 \begin{algorithm}

%\caption{Variable Step RAGD, \cite{ZS18}}
%\label{alg:C}
\begin{algorithmic}
\STATE {input : $\overrightarrow{\vec{x}}_0, \overrightarrow{\vec{v}}_0, 0<c \le\alpha_t< \frac{1}{L}, \beta_t >0$,$\delta >0$}, 
\\
\REPEAT
\STATE Compute $s_t\in(0,1)$ from the equation $s_t^2=\alpha_t((1-s_t)\gamma_t+s_t\mu)$.
\\
Set $\bar{\gamma}_{t+1}=(1-s_t)\gamma_t+s_t\mu$, $\gamma_{t+1}=\frac{1}{1+\beta_t}\bar{\gamma}_{t+1}$
\\
For $i\in[n]$:
\\
Set $\vec{y}_{i,t}=\Exp_{\vec{x}_{i,t}}\left(\frac{s_t\gamma_t}{\gamma_t+s_t\mu}\Log_{\vec{x}_{i,t}}(\vec{v}_{i,t})\right)$
\\
Set $\vec{x}_{i,t+1}=\Exp_{\vec{y}_{i,t}}(-\alpha_t\grad f(\vec{y}_{i,t}))$
\\
Set $\vec{v}_{i,t+1}=\Exp_{\vec{y}_{i,t}}\left(\frac{(1-s_t)\gamma_t}{\bar{\gamma}_t}\Log_{\vec{y}_{i,t}}(\vec{v}_{i,t})-\frac{s_t}{\bar{\gamma}_t}\grad f(\vec{y}_{i,t})\right)$

\UNTIL{$\grad f(\vec{y}_{i,t}) \le \delta$} for all $i\in[n]$.
\end{algorithmic}
\caption{ RAGD}
\end{algorithm}
\end{comment}

%\subsection{Saddle Avoidance Analysis}
Note that the update rule of RAGD can be written as the composition of three maps, i.e., denote 
\begin{align}
\vec{y}(\vec{x},\vec{v})&=\Exp_{\vec{x}}\left(\frac{s_t\gamma_t}{\gamma_t+s_t\mu}\Log_{\vec{x}}(\vec{v})\right),
\\
F(\vec{y})&=\Exp_{\vec{y}}(-\alpha_t\grad f(\vec{y})),
\\ 
 G(\vec{y},\vec{v})&=\Exp_{\vec{y}}\left(\frac{(1-s_t)\gamma_t}{\bar{\gamma}_t}\Log_{\vec{y}}(\vec{v})-\frac{s_t}{\bar{\gamma}_t}\grad f(\vec{y})\right).
 \end{align}
  Therefore, the update rule, which can be viewed as a map $\psi(\vec{x},\vec{v}):M\times M\rightarrow M\times M$, can be written in the following way:
\begin{align}\label{augmented map}
\psi(\vec{x},\vec{v})\overset{\text{def}}{=}\left(F(\vec{y}),G(\vec{y},\vec{v})\right)
=\left(F(\vec{y}(\vec{x},\vec{v})),G(\vec{y}(\vec{x},\vec{v}),\vec{v})\right).
\end{align}
If $\vec{x}^*=\vec{v}^*$ where $\vec{x}^*$ is a critical point, then the intermediate variable $\vec{y}(\vec{x}^*,\vec{v}^*)=\vec{x}^*$, and $\psi(\vec{x}^*,\vec{v}^*)=(\vec{x}^*,\vec{v}^*)$. This means that if $\vec{x}^*$ is a critical point of $f:M\rightarrow\mathbb{R}$, the point $(\vec{x}^*,\vec{x}^*)$ is a stationary point of the dynamical system defined by iterations of $\psi$ on $M\times M$.

The main technical step is the proof of "Stable-manifold theorem" for the dynamical system defined by iteration of the augmented map (\ref{augmented map}), which is a non-autonomous discrete-time dynamical system. We state the new stable manifold theorem as follows:
\begin{theorem}\label{generalSMT}
For a fixed dimension $d\ge 2$, there exists $s\in[d]$, and positive constant $K_1,K_2<1$, such that for all $t\in\mathbb{N}$, the following holds
\begin{align}
\sup_{1\le j\le s}\{\mathcal{L}_j(t)\}\le K_1, \ \ \ \ \ 
\sup_{s+1\le j\le d}\{\mathcal{L}_j^{-1}(t)\}\le K_2.
\end{align}
Suppose $\vect{\eta}(t,\vec{x})$ satisfies that $\vect{\eta}(t,\vec{0})=\vec{0}$ and for each $\epsilon>0$, there exists a neighborhod of $\vec{0}$ such that
\[
\norm{\vect{\eta}(t,\vec{x})-\vect{\eta}(t,\vec{y})}\le\epsilon\norm{\vec{x}-\vec{y}}.
\]
Then the dynamical system,
\[
\vec{x}_{t+1}=\vec{L}\vec{x}_t+\vect{\eta}(t,\vec{x}_t),
\]
where $\vec{L}=\diag\{\mathcal{L}_j(t)\}$,
has a local stable manifold in a small neighborhood of $\vec{0}$.
\end{theorem}

Note that the Center-Stable Manifold Theorem only secures the existence of local stable manifold, which implies that there exists a small neighborhood of an unstable fixed point $\vec{x}^*$ such that the all initial conditions taken from this neighborhood and converges to $\vec{x}^*$ must belong to a lower dimensional curved space, thus has measure zero. The local result can be extended to global through an argument in Appendix, for the case when the dynamical system is non-autonomous. In particular, if the update rule $\psi$ of the algorithm is not time-dependent such as the cases of constant step, the above result holds trivially by letting $\psi(t,\vec{x})\equiv\psi(\vec{x})$. We give the following general version on manifold also for further referrings.

\begin{proposition}\label{globalmeasurezero}
Let $\psi(t,\vec{x}):\mathbb{N}\times M\rightarrow M$ be a dynamical system on a finite dimensional manifold $M$. Suppose that $\psi(t,\cdot)$ is a local diffeomorphism on $M$ for each $t\in \BN$. Let $\mathcal{A}^*$ be the set of unstable fixed points. Then the global stable set $W^s(\mathcal{A}^*)$ of $\mathcal{A}^*$ has measure of zero with respect to the volume measure induced by the Riemannian metric on $M$.
\end{proposition} 
Here the stable and unstable set of a point $\vec{x}$ are defined to be
\begin{itemize}
\item Stable set of $\vec{x}$: \[W^s(\vec{x})=\{\vec{z}\in M:\lim_{t\rightarrow\infty} d(\psi(t,\vec{z}),\psi(t,\vec{x}))\rightarrow0\}\]
\item Unstable set of $\vec{x}$:
\[
W^u(\vec{x})=\{\vec{z}\in M:\lim_{t\rightarrow-\infty} d(\psi(t,\vec{z}),\psi(t,\vec{x}))\rightarrow 0\}.
\]
\end{itemize}

To prove following theorem \ref{SMT:RAGD2}, we leverage the property of the dynamical system induced by the map $\psi$ constructed in (\ref{augmented map}). Its orbit is the same as the orbit of RAGD and it also satisfies all other conditions in the center-stable manifold theorem. As mentioned above, the theorem tells us in a local neighborhood of a strict saddle point $\vec{x}^{\ast}$, all points in this neighborhood that will converge to $\vec{x}^{\ast}$ under $\psi$ lie in a submanifold, and it's dimension equals to the number of the eigenvalues of the Jacobian  of  $\psi$ at $\vec{x}^{\ast}$ that are less than $1$. Then we verify $DF(\vec{x}^{\ast},\vec{x}^{\ast})$ has a eigenvalue greater or equal to $1$. So all initial points in the local neighborhood  that will converge to $\vec{x}^{\ast}$ lies on a low dimension space, thus has measure $0$. Now we assume critical points are uncountable.  Since points converge to a strict saddle point under RAGD will finally fall into a local neighborhood of this strict saddle point, the set consisted by points in the whole space that will converge to a strict saddle point is a countable union of measure $0$ sets, which is also measure $0$. If critical points are uncountable, the Lindelof's lemma, which says every open cover there is a countable subcover, leads us to the fact that the set of initial conditions that converge to the set of saddle points (uncountable) is a countable union of measure zero sets, thus, a measure zero set. Formally, it is stated as follows.

\begin{theorem}\label{SMT:RAGD2}
Suppose $f:M\rightarrow\mathbb{R}$ is geodesically $\mu$-convex in neighborhood of local minima. Then the set of initial points $(\vec{x}_0,\vec{v}_0)\in M\times M$ converging to saddle points of $f$ under the Riemannian Accelerated Gradient Descent \cite{ZS18} has measure (induced by the product metric) of zero.
\end{theorem}

The proof is left in Appendix.

\begin{corollary}\label{distributed}
Let $\mathcal{M}=M_1\times...\times M_n$. Suppose that $f:\mathcal{M}\rightarrow\mathbb{R}$ is geodesically $\mu$-convex in neighborhood of local minima. Then the set of initial points $(\overrightarrow{\vec{x}}_0,\overrightarrow{\vec{v}}_0)\in\mathcal{M}\times\mathcal{M}$ converging to saddle points of $f$ has measure (induced by the product metric) of zero.
\end{corollary}
The proof is completed by trivially considering $\mathcal{M}$ as the manifold in Theorem \ref{SMT:RAGD2}, all the structures of the proof can be carried over. The following corollary is immediate from Theorem \ref{SMT:RAGD2} and Corollary \ref{distributed}.
%The theoretical guarantee in saddle avoidance of Riemannian Accelerated Gradient Descent implies that the Accelerated Multiplicative Weights Update avoids saddle points almost always, i.e., 

%The following corollary is immediate from Theorem \ref{SMT:RAGD2} and Corollary \ref{distributed}.

\begin{corollary}\label{AMWU}
Let $\mathcal{M}=\Delta_1\times...\times\Delta_n$, where 
\[
\Delta_i=\left\{\vec{x}_i\in\mathbb{R}^{d_i}:\sum_{s=1}^{d_i}x_{is}=1,x_{is}>0\right\}.
\]
Suppose $f:\mathcal{M}\rightarrow\mathbb{R}$ is $C^2$ and geodesically convex w.r.t. the product Shahshahani metric. Then the Accelerated Multiplicative Weights Update algorithm avoids interior saddle points almost always, i.e., randomly choose an initial point, the algorithm will avoid interior saddle points with probability one. 
\end{corollary}

\begin{remark}
In many applications, the constraint is compact, i.e., $\Delta_i=\left\{\vec{x}_i\in\mathbb{R}^{d_i}:\sum_{s=1}^{d_i}x_{is}=1,x_{is}\ge0\right\}$, and we might be interested in convergence to second-order stationary points \cite{ICMLPPW19,Meisam2}. With Theorem \ref{generalSMT}, Lemma 3.1 and Lemma 3.2 in \cite{ICMLPPW19}, one can extend the result of Corollary \ref{AMWU} to convergence to second-order stationary points with compact constraints. 
\end{remark}

%% file: expnew.tex
\section{Experiments}
\paragraph{Algorithms for comparison.} The experiments are designed to understand the behavior of A-MWU compared to classic MWU and Accelerated Mirror Descent (A-MD) proposed by \cite{KBB}, especially their behaviors near a saddle point. Since the performance of A-MD is controlled by the parameter $r$ in \cite{KBB}, i.e., the larger $r$ ends up with a smoother curve of convergence, but with a slower rate of convergence. In order to understand exactly how A-MWU performs, we choose different $r$ in Accelerated Mirror Descent while keep the step-size of A-MWU and A-MD the same. 

\begin{figure}[htp]
\centering
\subfigure[A-MWU, MWU, A-MD, step-size=0.001 ]{
\includegraphics[clip,width=0.7\columnwidth]{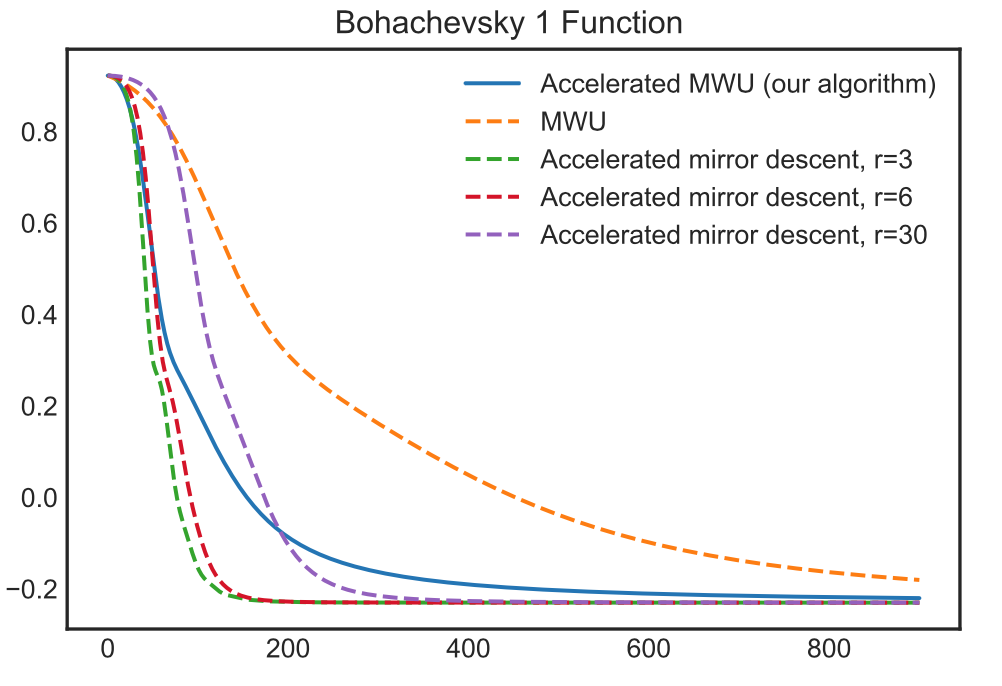}
}

\subfigure[Trajectories]{
\includegraphics[clip,width=0.7\columnwidth]{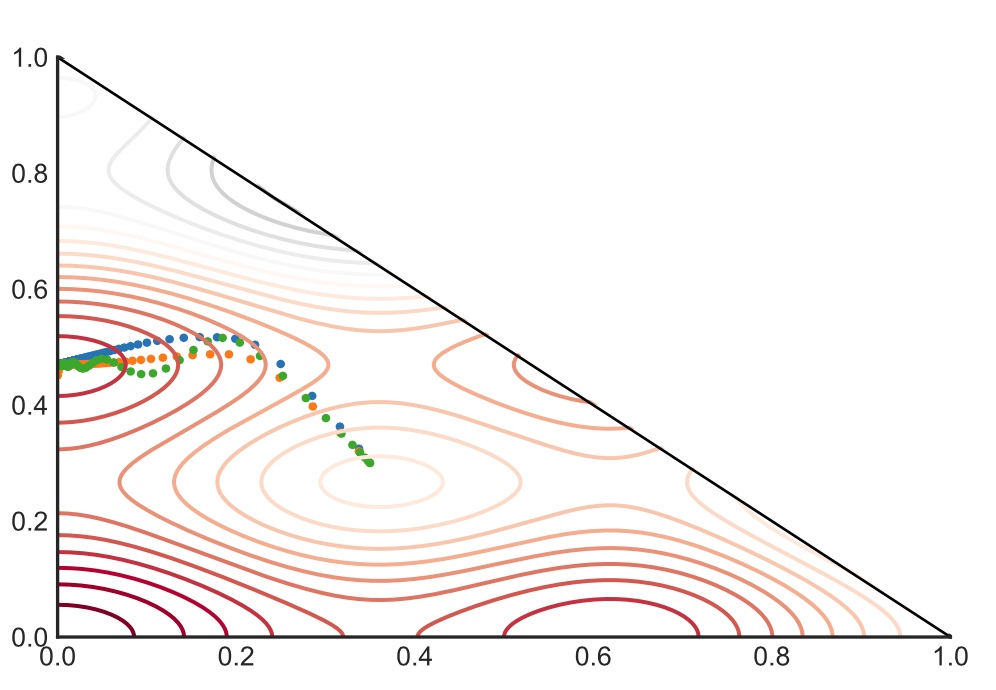}
}

\caption{Bohachevsky function}
\label{Boha}
\end{figure}

\begin{figure}[htp]
\centering
\subfigure[A-MWU, MWU, A-MD, step-size=0.01]{
  \includegraphics[clip,width=0.7\columnwidth]{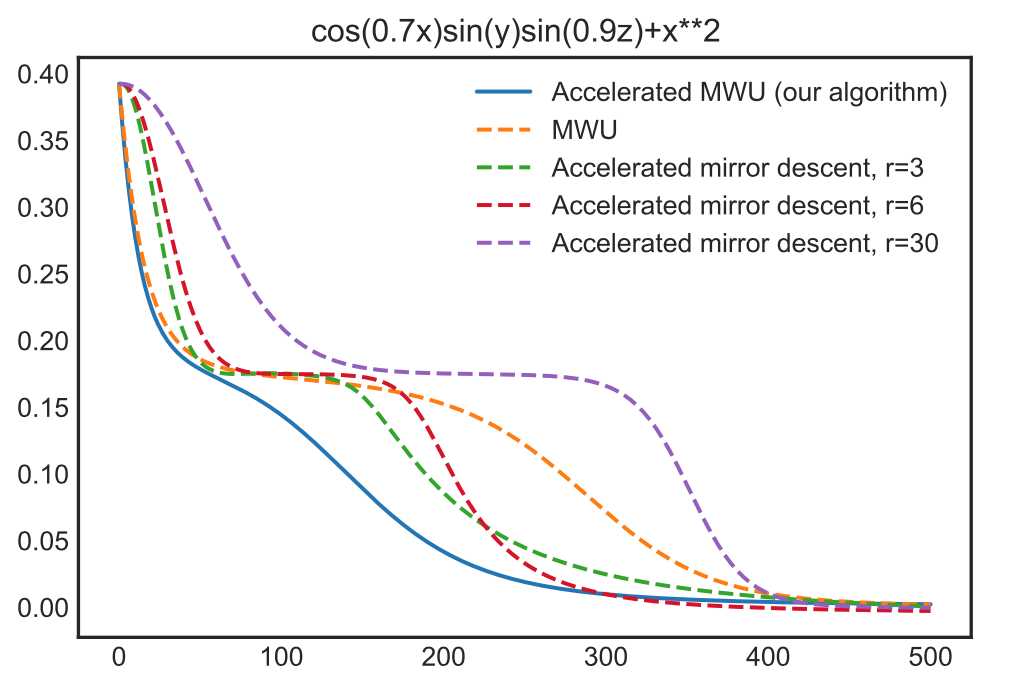}
}

\subfigure[Trajectories]{
  \includegraphics[clip,width=0.7\columnwidth]{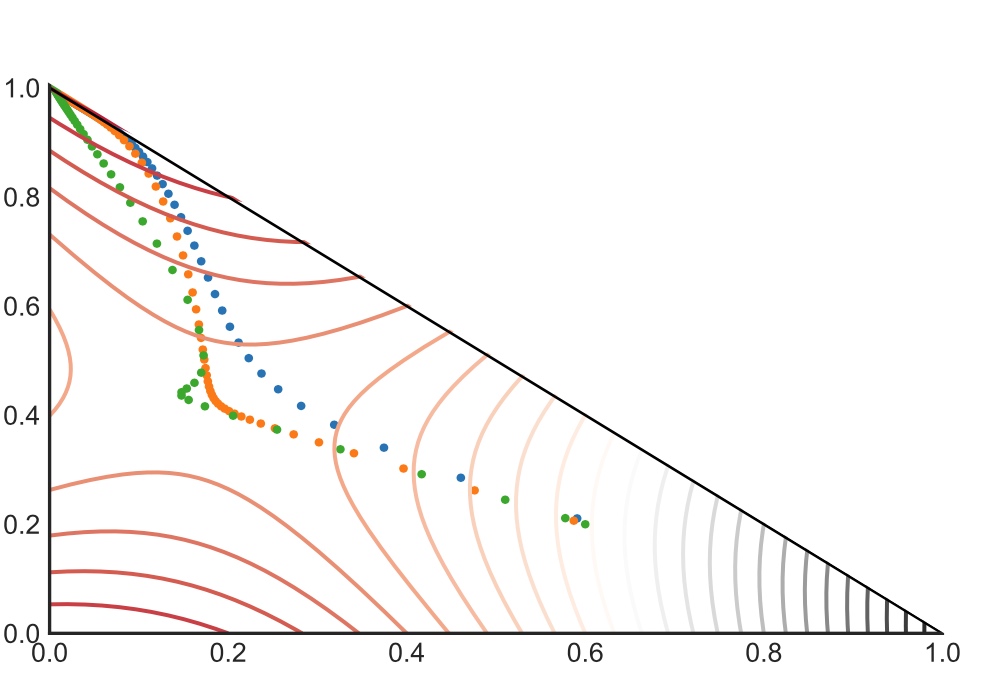}
}

\caption{Test function 2:$\cos(0.7x)\sin(y)\sin(0.9z)+x^2$}
\label{trig2}
\end{figure}

\begin{figure}[htp]
\centering
\subfigure[A-MWU, MWU, step-size = 0.001]{
\includegraphics[clip,width=0.7\columnwidth]{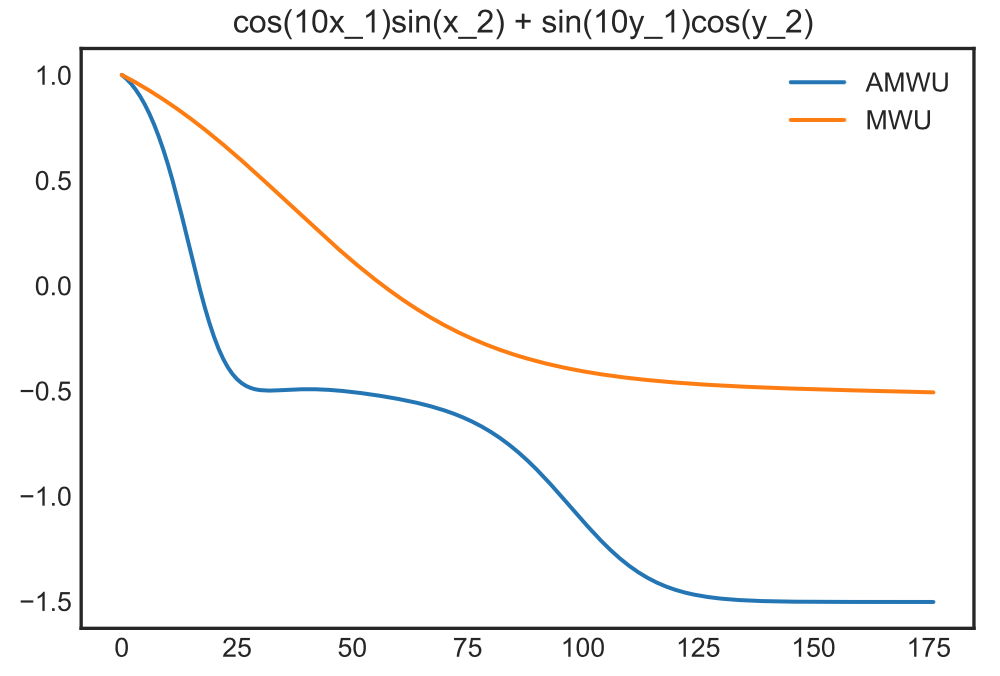}
}
\subfigure[Trajectories]{
\includegraphics[clip,width=0.7\columnwidth]{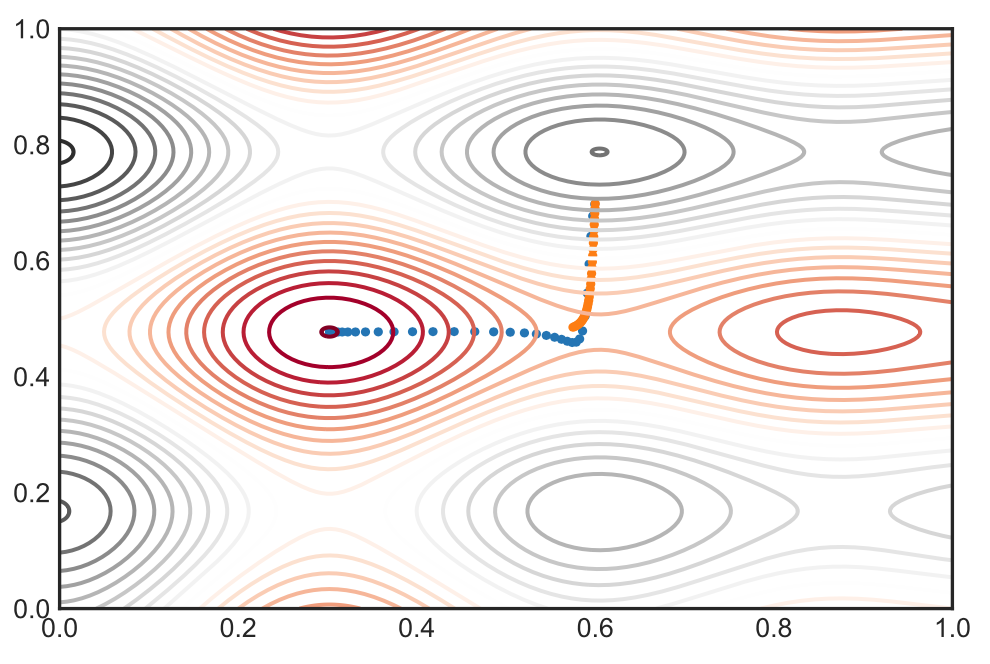}
}
\caption{Product of 1-simplices.}
\label{twosimplex}
\end{figure}

\paragraph{Parameter setting.} In our experiments, we set the parameters as follows:
\begin{itemize}
\item $\beta=0.001$ and $\mu=1$ for Rosenbrock function, Figure \ref{Rosenb}.
\item $\beta=0.1$ and $\mu=1$ for Bohachevsky function, Figure \ref{Boha}.
\item $\beta=0.1$ and $\mu=0.2$ for Test function 1, Figure \ref{trig1}.
\item $\beta=0.001$ and $\mu=0.001$ for Test function 2, Figure \ref{trig2}.
\item $\beta=0.1$ and $\mu=0.5$ for two-agent, Figure \ref{twosimplex}.
\end{itemize}
In general, by definition of geodesic convexity, the principle of choosing $\mu$ is to estimate a lower bound of the actual geodesic convexity parameter $\mu$. Since when $\mu>\mu_0$, then the convexity inequality holds for $\mu_0$ if it holds for $\mu$. To achieve fast local convergence rate, one can choose $\beta\approx\frac{1}{5}\sqrt{\mu/L}$ according to \cite{ZS18}. There is a further discussion on parameters in Appendix.

\paragraph{Efficiency in convergence.} For the local convergence behavior, we use Rosenbrock function and Bohachevsky function as the test functions. 
\begin{itemize}
\item The Rosenbrock function: $(0.5 - x)^2 + 0.25 (y-x^2)^2 + x+y+z-1$, with initial points $(0.2,0.4,0.4)$, and the step sizes are $0.01$ for A-MWU, MWU and  A-MD, see Figure \ref{Rosenb}.

\begin{figure}[htp]
\centering
\subfigure[A-MWU, MWU, A-MD, step-size=0.01]{
  \includegraphics[clip,width=0.7\columnwidth]{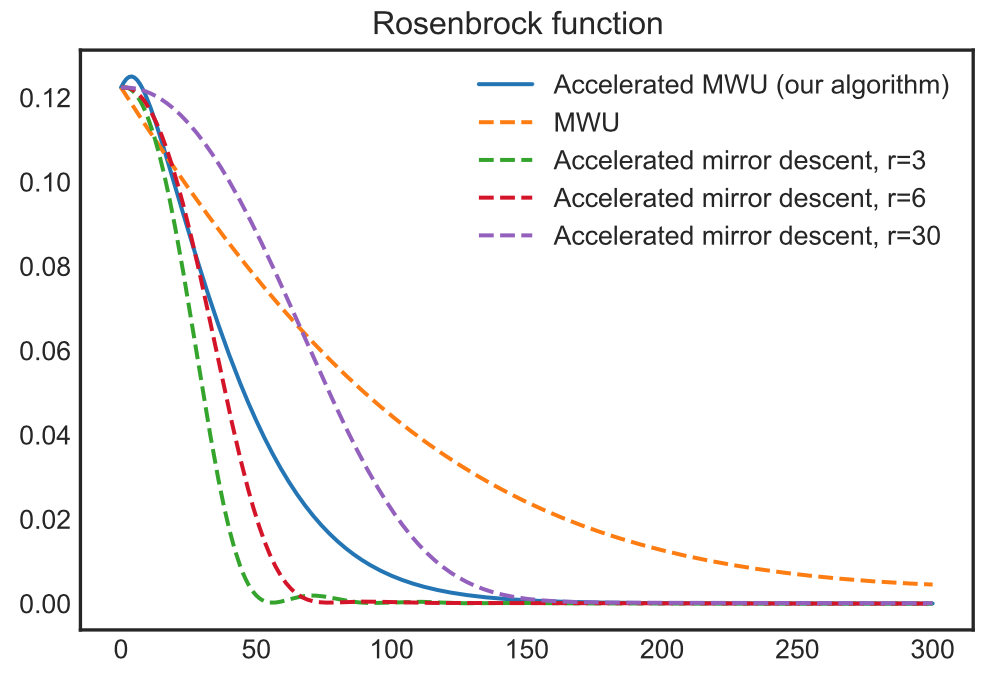}
}

\subfigure[Trajectories]{
  \includegraphics[clip,width=0.7\columnwidth]{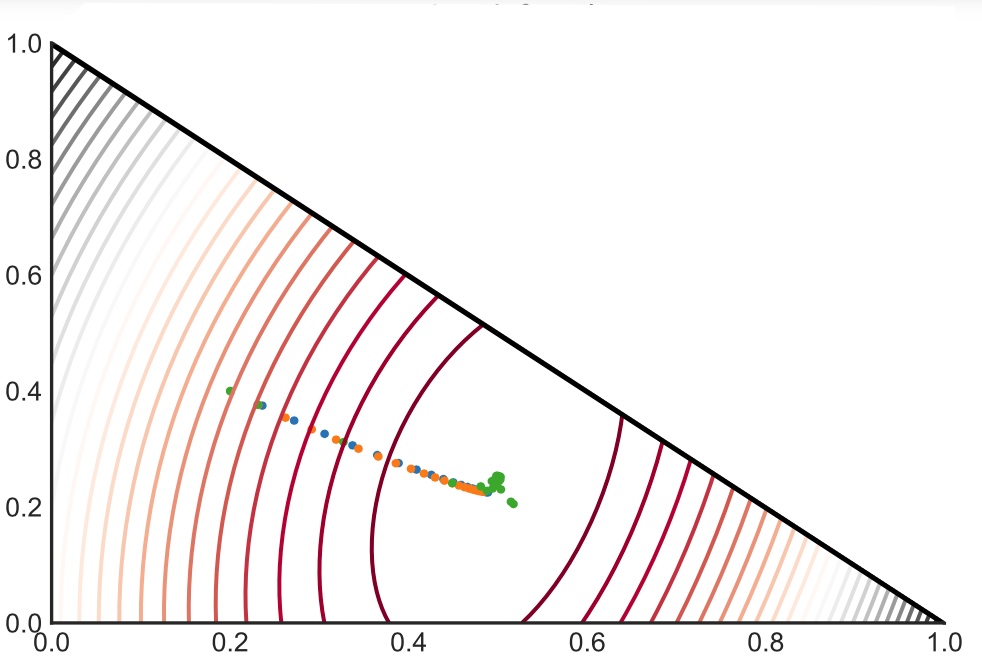}
}

\caption{Rosenbrock function.}
\label{Rosenb}
\end{figure}

\item The Bohachevsky function: $x^2+2y^2-0.3\cos(3\pi x)-0.4\cos(4\pi y)+x+y+z-1$, with initial point $(0.35,0.3,0.35)$, and step-size 0.001, see Figure \ref{Boha}.
\end{itemize}

Compared to the usual definitions of $2$-dimension Rosenbrock function and Bohachevsky function, our expression contains an additional term $x+y+z-1$ to make the function well defined in $\mathbb{R}^3$. But in simplices, $x+y+z-1 = 0$ always holds, thus our definition boils down to the original definitions.

\begin{comment}
\begin{figure}[htp]
\centering
\subfigure[A-MWU, MWU, A-MD, step-size=0.01]{
  \includegraphics[clip,width=0.7\columnwidth]{rosenbrock_1.png}
}

\subfigure[Trajectories]{
  \includegraphics[clip,width=0.7\columnwidth]{rosenbrock_2.jpg}
}

\caption{Rosenbrock function.}
\label{Rosenb}
\end{figure}
\end{comment}

\paragraph{Efficiency in escaping saddle points.}
We compare the curves of convergence and trajectories of A-MWU, A-MD and MWU for the function with many saddle points, in order to illustrate the different efficiency of three algorithms in escaping saddle points.

\begin{itemize}
\item Test function 1: $\cos(8.5x)\sin(8.5(y-0.4))+\sin(8.5z)$, with initial point $(0.42,0.24,0.33)$, and the step-size is $0.005$.

\begin{figure}[htp]
\centering
\subfigure[A-MWU, MWU, A-MD, step-size=0.005]{
  \includegraphics[clip,width=0.7\columnwidth]{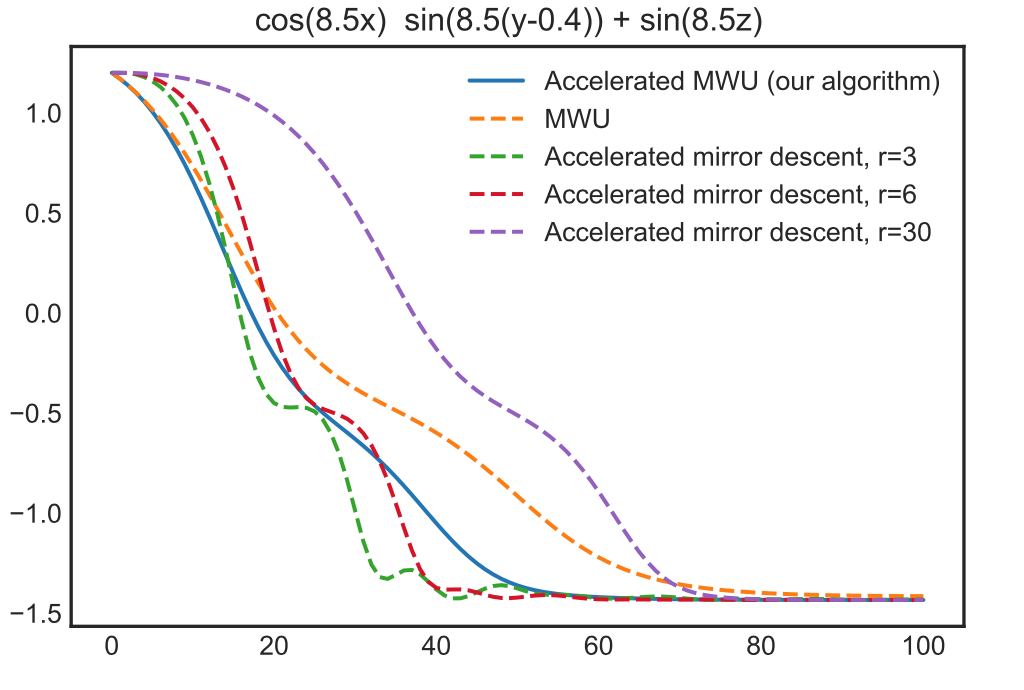}
}

\subfigure[Trajectories]{
  \includegraphics[clip,width=0.7\columnwidth]{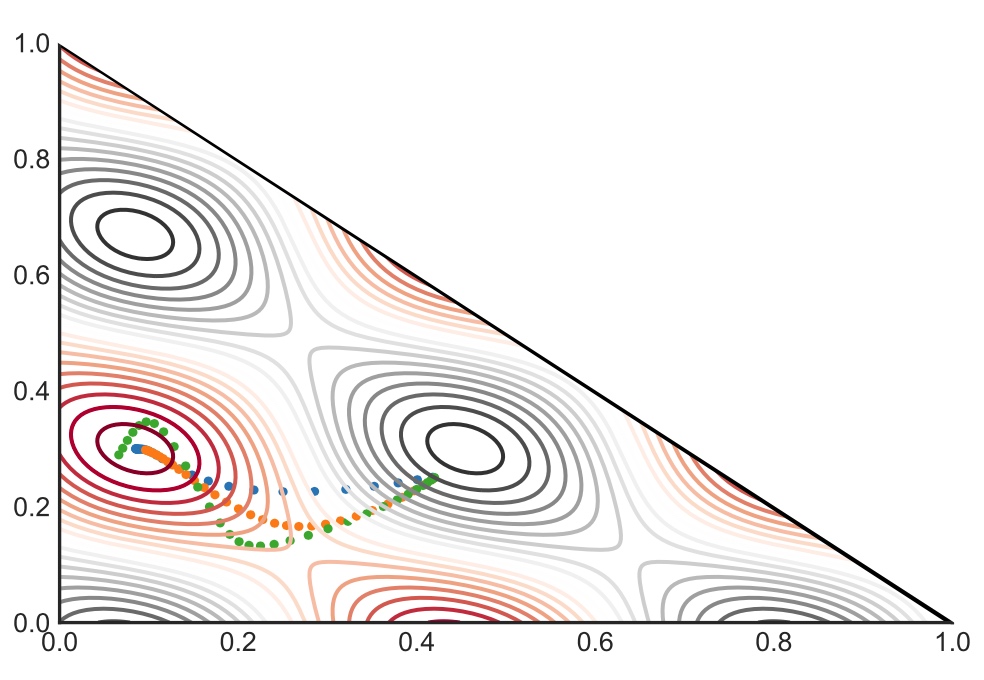}
}

\caption{Test function 1:$\cos(8.5x)\sin(8.5(y-0.4))+\sin(8.5z)$}
\label{trig1}
\end{figure}

\item Test function 2: $\cos(0.7x)\sin(y)\sin(0.9z)+x^2$, with initial point $(0.6,0.2,0.2)$, and the step-size is $0.05$
\end{itemize}

\begin{comment}
\begin{figure}[htp]
\centering
\subfigure[A-MWU, MWU, A-MD, step-size=0.005]{
  \includegraphics[clip,width=0.7\columnwidth]{triangle1_1.png}
}

\subfigure[Trajectories]{
  \includegraphics[clip,width=0.7\columnwidth]{triangle1_2.jpg}
}

\caption{Test function 1:}
\label{trig1}
\end{figure}
\end{comment}

\paragraph{Multiple simplices.} We give a final illustration on the behavior of A-MWU in two-agent case. Suppose $(x_1,x_2)\in[0,1]\times[0,1]$ and satisfies $x_1+x_2=1$, $(y_1,y_2)\in[0,1]\times[0,1]$ and $y_1+y_2=1$. The function is of the form $\cos(10x_1)\sin(x_2)+\sin(10y_1)\cos(y_2)$. 

\paragraph{Experimental conclusions.} The experimental results indicate the following:
\begin{itemize}
\item As expected, all experiments have verified that A-MWU has a better convergence behavior and efficiency in escaping saddle points compared to classic MWU.
\item Compared to Accelerated Mirror Descent with entropy regularizer \cite{KBB}, the proposed algorithm A-MWU has the similar acceleration effect as A-MD with small $r$ value, and outperforms A-MD with larger $r$ value. See Figure \ref{Rosenb} and \ref{Boha}.
\item From Figure \ref{trig1}-(a) and Figure \ref{trig2}-(a), we believe that A-MWU has potential to outperform A-MD significantly in escaping saddle points.   
\end{itemize}